\documentclass[review,onefignum,onetabnum]{siamart171218}

\usepackage{lipsum}
\usepackage{amsfonts}
\usepackage{graphicx}
\usepackage{epstopdf}

\ifpdf
  \DeclareGraphicsExtensions{.eps,.pdf,.png,.jpg}
\else
  \DeclareGraphicsExtensions{.eps}
\fi

\newsiamremark{remark}{Remark}
\newsiamremark{hypothesis}{Hypothesis}
\crefname{hypothesis}{Hypothesis}{Hypotheses}
\newsiamthm{claim}{Claim}

\title{CLINN: Conservation Law Informed Neural Network for Approximating Discontinuous Solutions}

\author{Weiheng Zeng\thanks{LMIB and School of Mathematical Sciences, Beihang University (\email{harry\_zeng@buaa.edu.cn})}
\and Kun Wang\thanks{LMIB and School of Mathematical Sciences, Beihang University (\email{wangkun@buaa.edu.cn})}
\and Ruoxi Lu\thanks{LMIB and School of Mathematical Sciences, Beihang University (\email{sorenblu@buaa.edu.cn})}
\and Tiegang Liu\;\textsuperscript{\faEnvelopeO}\;\thanks{LMIB and School of Mathematical Sciences, Beihang University (\email{liutg@buaa.edu.cn})}}

\makeatletter
\newcommand*{\addFileDependency}[1]{
  \typeout{(#1)}
  \@addtofilelist{#1}
  \IfFileExists{#1}{}{\typeout{No file #1.}}
}
\makeatother

\ifpdf
\hypersetup{
  pdftitle={CLINN: Conservation Law Informed Neural Network for Approximating Discontinuous Solutions},
  pdfauthor={W. Zeng, K. Wang, R. Lu, T. Liu}
}
\fi

\begin{document}

\maketitle

\begin{abstract}
Physics-informed Neural Network (PINN) faces significant challenges when approximating solutions to conservation laws, particularly in ensuring conservation and accurately resolving discontinuities. To address these limitations, we propose Conservation Law-informed Neural Network (CLINN), a novel framework that incorporates the boundedness constraint, implicit solution form, and Rankine-Hugoniot condition of scalar conservation laws into the loss function, thereby enforcing exact conservation properties. Furthermore, we integrate a residual-based adaptive refinement (RAR) strategy to dynamically prioritize training near discontinuities, substantially improving the network's ability to capture sharp gradients. Numerical experiments are conducted on benchmark problems, including the inviscid Burgers equation, the Lighthill-Whitham-Richards (LWR) traffic flow model, and the Buckley-Leverett problem. Results demonstrate that CLINN achieves superior accuracy in resolving solution profiles and discontinuity locations while reducing numeral oscillations. Compared to conventional PINN, CLINN yields a maximum reduction of 99.2\% in mean squared error (MSE). 
\end{abstract}

\begin{keywords}
  conservation law, artificial intelligence, physics-informed neural network, discontinuous solution
\end{keywords}


\section{Introduction}
\label{sec:intro}

The shock waves, generated by aircraft during supersonic flight in the aerospace field, can be mathematically characterized as discontinuous solutions to hyperbolic conservation laws. Such solutions pose considerable challenges in simulations due to their numerical instability. As a result, the development of robust numerical methods capable of solving hyperbolic conservation laws while accurately capturing and resolving shock discontinuities has long been an enduring challenge in the field of computational fluid dynamics (CFD).

High-resolution numerical methods, represented by the Weighted Essentially Non-Oscillatory (WENO) scheme \cite{jiang1996efficient, friedrich1998weighted, hu1999weighted} and the Discontinuous Galerkin (DG) method \cite{reed1973triangular, cockburn1989tvb, cockburn1991runge}, are based on function space approximation theory or Taylor expansions. These techniques approximate either the integral operator of the weak form or the differential operator of the strong form of partial differential equations (PDEs) by incorporating high-order residual terms, thereby achieving high numerical accuracy. Through shock-capturing strategies and specialized treatments for discontinuities, they mitigate—though do not fully eliminate—numerical instabilities near discontinuities. Nevertheless, such methods are not without limitations, including elevated computational costs and degradation of solution accuracy in the vicinity of discontinuities.

In recent years, the emerging paradigm of artificial intelligence for PDE (AI4PDE) has garnered significant attention in computational mathematics. Current AI4PDE methodologies can be broadly classified into two distinct categories: (I) hybrid approaches that leverage deep learning to enhance traditional high-resolution numerical schemes in computational simulations, and (II) end-to-end frameworks that employ deep neural networks to approximate PDE solutions directly.

For the first category of approaches, a prominent strategy involves developing AI-based troubled-cell indicators to enhance numerical simulations. These AI-driven indicators address a critical limitation of conventional approaches, which frequently misclassify extremal points and high-gradient regions as discontinuities. For instance, Ray et al. \cite{ray2018artificial} introduced a neural network-based detection method that utilizes solution values at cell boundaries and integral averages of adjacent cell solutions as input features. However, this approach suffered from computational inefficiency due to its large network architecture. To overcome this limitation, Feng et al. \cite{feng2020characteristic, feng2021characteristic} proposed a novel shock indicator employing artificial neurons (ANs) that capitalizes on the mathematical properties of conservation laws while achieving remarkable efficiency with a minimal network structure of only three neurons. When implemented as a replacement for the total variable bounded (TVB) indicator in DG post-processing, this AN-based indicator exhibited superior numerical performance.

For the second category of approaches, a predominant strategy involves training neural networks to satisfy the governing PDE constraints, thereby enabling direct approximation of the solution. The physics-informed neural network (PINN) framework, introduced by Raissi et al. \cite{raissi2019physics}, represents a seminal advancement in this direction. This approach incorporates the governing PDE equations along with their initial and boundary conditions into the neural network's loss function, effectively embedding physical principles into the deep learning architecture. However, when applied to conservation law equations, the PINN framework demonstrates limited effectiveness due to two fundamental characteristics distinct from conventional PDEs: (I) their solutions must strictly satisfy conservation properties, and (II) they typically admit weak solutions that may contain discontinuities. Although it has been proven that multi-layer perceptrons (MLPs) can approximate any Borel measurable function \cite{cybenko1989approximation, hornik1989multilayer}, including discontinuous solutions of conservation law equations, the effectiveness of neural networks in approximating such discontinuous functions is still unsatisfactory in practice.

To guarantee that neural network predictions satisfy  conservation properties, it is necessary to embed conservation-related formulations into the neural network architecture to guide the network's training. For instance, Zhang et al. \cite{zhang2022implicit} proposed the implicit form neural network (IFNN), which incorporates the implicit solutions of conservation laws into the loss function. Liu et al. \cite{liu2024discontinuity} developed PINN with equation weight (PINN-WE), integrating the Rankine-Hugoniot jump condition \cite{chang1989riemann} into the loss function. Pater et al. \cite{patel2022thermodynamically} introduced thermodynamically consistent PINN (TC-PINN), embedding both the total variation diminishing (TVD) property \cite{harten1997high} and Lax entropy conditions \cite{lax1973hyperbolic} into the loss function. 

To approximate discontinuous solutions, Lu et al. \cite{lu2021deepxde} developed DeepXDE, a deep learning library for solving PDEs, which employs a residual-based adaptive refinement (RAR) algorithm to detect and incorporate problematic points (particularly near discontinuities) into the training set. Ferrer et al. \cite{ferrer2024gradient} proposed gradient-annihilated PINN (GA-PINN), which mitigates the influence of high-gradient regions by adaptively reducing their loss weighting, thereby prioritizing optimization in smooth solution domains. Chen et al. \cite{chen2022bridging, chen2023random} introduced the random feature method (RFM), approximating PDE solutions through linear combinations of step functions and random feature functions. De Ryck et al. \cite{de2024wpinns} presented weak PINN (wPINN), replacing strong-form PDEs with their weak formulations in the loss function to better handle discontinuous solutions. Sun et al. \cite{sun2024lift} proposed the lift-and-embed learning method (LELM), which embeds low-dimensional discontinuous functions into high-dimensional continuous functions, thereby transforming the problem into solving for a continuous solution. Liu et al. \cite{liu2025deepoly} introduced the DeePoly method, which combines the global approximation capability of neural networks with the local approximation strength of polynomials to enhance the accuracy in approximating local discontinuities. See also \cite{abbasi2025challenges} for more PINN-based methods addressing non-smooth problems. 

Notably, existing approaches still demonstrate limited effectiveness when modeling wave system interactions, revealing neural networks' persistent shortcomings in maintaining conservation properties and approximating discontinuities. In response to these persistent challenges, we propose a conservation laws informed neural network (CLINN), with the following key contributions:

\begin{itemize}
  \item We incorporate both the fundamental physical principles common to PDEs and distinctive features of conservation laws (e.g., implicit solution, boundedness constraint, Rankine-Hugoniot jump condition) into the network architecture and loss function, thereby better guiding the neural network to learn solutions of  conservation laws.
  \item We employ the AN-based shock indicator to detect discontinuities, combined with an improved RAR method to accelerate convergence in regions near discontinuities.
  \item We conduct numerical experiments on the inviscid Burgers' equation, the LWR traffic flow model, and the non-convex/non-concave problem, including the Buckley-Leverett equation, covering complex wave interaction scenarios. Results demonstrate that CLINN achieves up to a reduction of 99.2\%  in Mean Squared Error (MSE) compared to PINN.
\end{itemize}

\section{Preliminaries}
\label{sec:pre}
In this section, we investigate the fundamental properties of conservation laws and their solutions, along with the underlying principles and architecture of PINN.

\subsection{Scalar conservation laws}
This paper focuses on the initial value problem of scalar hyperbolic conservation laws
\begin{equation}
\begin{cases}
\mathcal P(u(\boldsymbol x,t)):=\big(\partial_t u+\nabla_{\boldsymbol x} \cdot \boldsymbol f(u)\big)\Big|_{(\boldsymbol x,t)}=0, &\boldsymbol x\in\Omega\subseteq\R^d, t>0; \\ 
u(\boldsymbol x,0)=u_0(\boldsymbol x), & \boldsymbol x\in\Omega\subseteq\R^d;
\end{cases}
\label{pde}\end{equation}
where $u(\boldsymbol x,t)\in\R$ represents the conserved variable to be solved, and $\boldsymbol f:=(f_1, \cdots, f_d)$ $\in C^\infty(\R^d), d\in\mathbb{N}$ denotes the numerical flux vector. We additionally assume the initial value function $u_0(\boldsymbol x)$ is bounded on $\R$.

Based on Eq. (\ref{pde}) and characteristic theory, we summarize the following key conclusions to be referenced in later sections. See \cite{chang1989riemann} for detailed proof.

\begin{proposition}
Denoting $\boldsymbol \lambda(u):=\boldsymbol f'(u)$, then we have following properties.

\begin{itemize}
\item \textbf{Implicit Solution}: Eq.
    (\ref{pde}) admits the  solution in an implicit form
    \begin{equation}
        u=u_0(\boldsymbol x-\boldsymbol\lambda(u)t),\quad \forall (\boldsymbol x,t)\in D\subseteq\Omega\times(0,+\infty),
    \label{solution}\end{equation}
    where $D$ satisfies $u\in C^\infty(D)$. In other words, Eq.(\ref{solution}) holds only in regions where $u$ remains smooth.
    \item \textbf{Boundedness}: All solutions remain bounded, satisfying
    \begin{equation}
        \inf\limits_{x\in\Omega} u_0(\boldsymbol x)\leq u(\boldsymbol x,t)\leq \sup\limits_{x\in\Omega} u_0(\boldsymbol x),\quad\forall(\boldsymbol x,t)\in\Omega\times[0,+\infty).
    \label{bound}\end{equation}

    \item \textbf{Jump condition}: The characteristic lines of Eq. (\ref{pde}) exhibit linear trajectories. These lines may converge, leading to the formation of discontinuity surface $\Gamma:\ \Phi(\boldsymbol 
    x, t)=0$. The normal propagation velocity of the discontinuity surface $s(\boldsymbol x, t):=-\frac{\partial_t\Phi}{|\nabla_{\boldsymbol x }\Phi|}$ satisfies the Rankine-Hugoniot jump condition 
    \begin{equation}
    s=\boldsymbol n\cdot \frac{[\boldsymbol{f}(u)]}{[u]},
    \label{RH}
    \end{equation}
    where $\boldsymbol n:=\frac{\nabla_{\boldsymbol x }\Phi}{|\nabla_{\boldsymbol x }\Phi|}$ denotes the unit outer normal vector of $\Phi(\cdot, t)$, and $[\cdot]$ denotes the jump value of conserved variables or fluxes across the discontinuity. Denote
    \begin{equation}
        u^{\pm}(\boldsymbol x, t):=\lim\limits_{\varepsilon\to 0} u(\boldsymbol x\pm\varepsilon \boldsymbol n, t),
    \end{equation} then $[u]:=u^+-u^-, [\boldsymbol f(u)]:=\boldsymbol f(u^+)-\boldsymbol f(u^-)$.    
    
    \item \textbf{Entropy condition}: The physical solution is uniquely determined by both the weak form of Eq.(\ref{pde}) and the Oleinik\cite{oleinik1963discontinuous} entropy condition
    \begin{equation}
        \boldsymbol\lambda(u^+)\cdot \boldsymbol n\leq s\leq \boldsymbol\lambda(u^-)\cdot \boldsymbol n.
    \label{Oleinik}\end{equation}

\end{itemize}
\label{prop}
\end{proposition}

Since the solutions of Eq. (\ref{pde}) may develop discontinuities, it is meaningful to study the solutions of one-dimensional scalar conservation laws with piecewise constant initial data 
\begin{equation}u_0(x)=\begin{cases}u_l, & x\leq x_0;\\u_r, & x>x_0.\end{cases}\quad(u_l\neq u_r),\label{Riemann}\end{equation}
which is known as the Riemann problem. 

\begin{figure}
    \centering
    \includegraphics[width=0.85\linewidth]{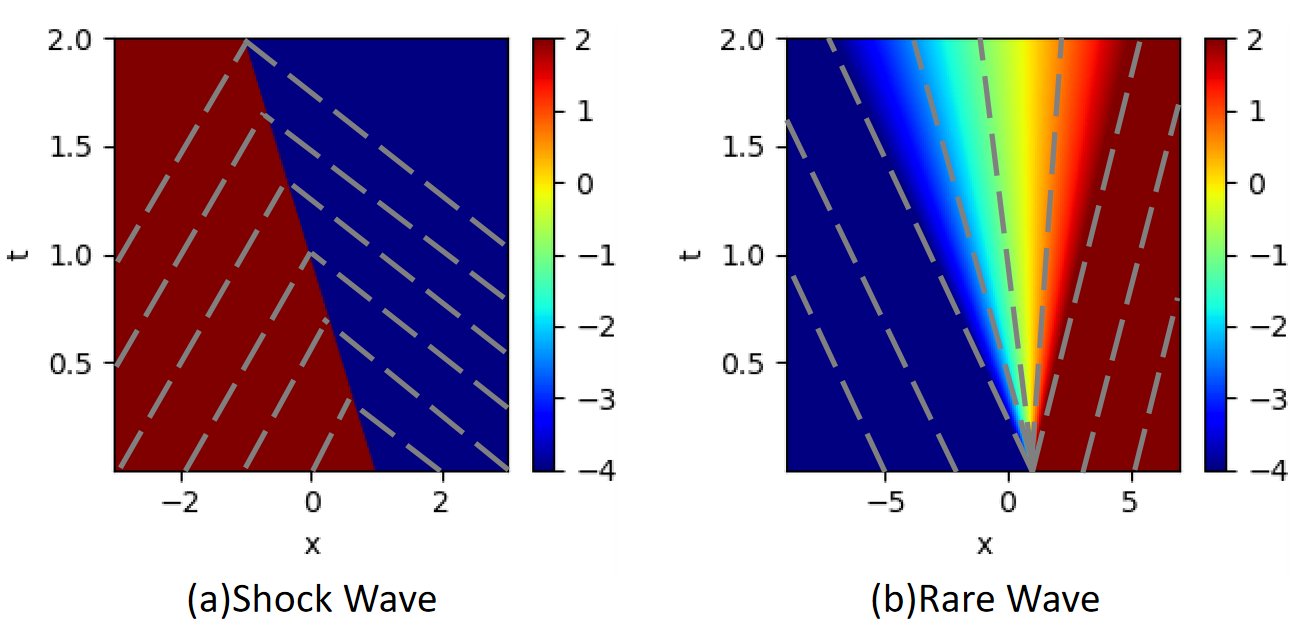}
    \caption{Solutions of Riemann problems with characteristic lines indicated by dashed lines. (a) Case containing a shock wave; (b) Case containing a centered rarefaction wave. Parameters are set to be: $f(u)=u^2/2,\ x_0=1,\ \min\limits_x u_0(x)=-4,\ \max\limits_x u_0(x)=2.$}
    \label{SW_RW}
\end{figure}

\begin{proposition}
    Supposing $f(u)$ is convex, i.e. $\lambda'(u)>0$. 
\begin{itemize}
    \item if $\lambda(u_l)>\lambda(u_r)$, then the characteristic lines originating from both constant-state regions converge, forming a shock wave at the intersection points, see Fig. \ref{SW_RW}(a). The solution is \begin{equation}u(x,t)=\begin{cases}u_l, x<x_0+st;\\ u_r, x>x_0+st.\end{cases}\label{S-case}\end{equation}
    where $s=\frac{f(u_l)-f(u_r)}{u_l-u_r}$, given by (\ref{RH}).
    
    \item if $\lambda(u_l)<\lambda(u_r)$, then the characteristic lines originating from both constant-state regions diverge, resulting in a centered rarefaction wave in the intermediate region, see Fig. \ref{SW_RW}(b). The solution is\begin{equation}u(x,t)=\begin{cases}u_l&x<x_0+\lambda(u_l)t;\\\lambda^{-1}\left(\dfrac{x-x_0}{t}\right)&x_0+\lambda(u_l)t\leq x\leq \lambda(u_r)t+x_0;\\u_r&x>x_0+\lambda(u_r)t.\end{cases}\label{R-case}\end{equation}
    where $\lambda^{-1}(\cdot)$ is the reverse function of $\lambda(\cdot)$.
\end{itemize}
\end{proposition}

The discontinuity is classified as a shock wave if neither side of the inequality in Eq. (\ref{Oleinik}) attains equality; otherwise, it is termed a contact discontinuity, which may occur when the flux is neither convex nor concave. Shock waves, contact discontinuities, and centered rarefaction waves constitute the fundamental wave types in hyperbolic conservation laws.

\subsection{Physics-informed neural network (PINN)}

The PINN framework \cite{raissi2019physics} incorporates the governing PDE and initial/boundary conditions into the loss function of the neural network, thereby not only imbuing the approximation with rigorous mathematical-physical constraints but also ensuring the predictive accuracy.

\begin{figure}[H] 
    \centering 
    \includegraphics[width=\textwidth]{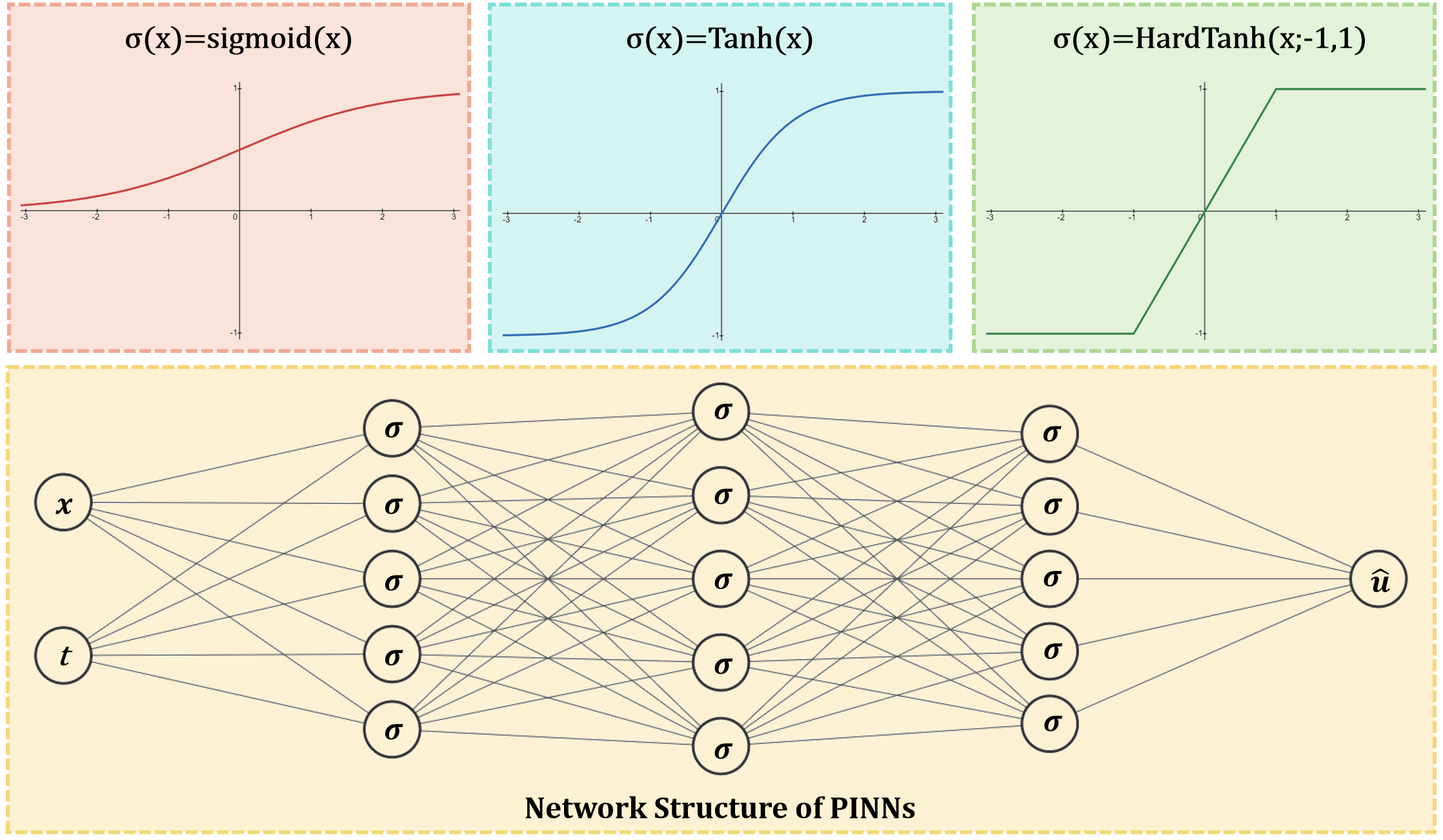} 
    \caption{The architecture and activation functions in PINN.} 
    \label{mlp} 
\end{figure}

The PINN employs a feed-forward network (FFN) structure, as is shown in Fig. \ref{mlp}. Consider a $K$-layer network with $n_k$ neurons in the $k$-th layer, the structure can be mathematically formulated as:
\begin{equation}
\begin{aligned}
    \widehat{u}^0(\boldsymbol x,t)&=[\boldsymbol x\quad t]^T\in\R^{n_0},n_0=d+1 \\
    \widehat{u}^k(\boldsymbol x,t)&=\sigma(\boldsymbol{W}^k \widehat{u}^{k-1}(\boldsymbol x,t) + \boldsymbol{b}^k)\in\R^{n_k}, k\in[1, K-1] \\
    \widehat{u}^{K}(\boldsymbol x,t)&=\boldsymbol{W}^{K} \widehat{u}^{K-1}(\boldsymbol x,t) + \boldsymbol{b}^{K}\in\R^{n_K}, n_K=1
\end{aligned}
\end{equation}
where $\boldsymbol{W}^k \in \mathbb{R}^{n_{k}\times n_{k-1}}$ denotes the weight matrix, $\boldsymbol{b}^k \in \mathbb{R}^{n_k}$ represents the bias vector, and $\sigma(\cdot)$ is the activation function that introduces nonlinearity, enabling the network to approximate complex nonlinear systems. The activation functions considered in this study include
\begin{equation}
\begin{aligned}
\mathrm{Tanh}(x)&=\dfrac{e^x-e^{-x}}{e^x+e^{-x}};\\
\quad\mathrm{Sigmoid}(x)&=\dfrac{1}{1+e^{-x}};\\ \mathrm{HardTanh}(x; c_1,c_2)&=\begin{cases}
c_1, & x<c_1\\ x, & c_1\leq x\leq c_2\\ c_2, &x>c_2
\end{cases} \ (c_1<c_2),
\label{act_func}
\end{aligned}
\end{equation}
which are shown in Fig. \ref{mlp}. We formally express this neural network as $\widehat{u}_{\sita}(x,t)$, where $\sita$ denotes the learnable parameter vector, comprising all elements of the weight matrices and bias vectors.

We now proceed to optimize the neural network parameters by minimizing an objective function, commonly referred to as the loss function. Our goal is to ensure that the neural network approximation $\widehat{u}_{\sita}(\boldsymbol x,t)$ of $u(\boldsymbol x,t)$ satisfies both the governing equation and initial/boundary conditions as closely as possible. Accordingly, for conservation laws, we define the composite loss function as
\begin{equation}
\begin{aligned}
\loss_{\mathrm{PINN}}(\boldsymbol x,t;\sita)&:=w_{\mathrm{GOV}}\loss_{\mathrm{GOV}}(\boldsymbol x,t;\sita)+w_{\mathrm{IC}}\loss_{\mathrm{IC}}(\boldsymbol x,t;\sita)+w_{\mathrm{IC}}\loss_{\mathrm{BC}}(\boldsymbol x, t,\sita);\\
\loss_{\mathrm{GOV}}(\boldsymbol x,t;\sita)&:=\sum\limits_{j: (\boldsymbol x_j,t_j)\in P_N} \loss_{\mathrm{GOV},j},\ \loss_{\mathrm{GOV},j}:=\dfrac{1}{|P_{N}|}\left|\mathcal P(\widehat u_{\sita, j})\right|^2;\\ 
\loss_{\mathrm{IC}}(\boldsymbol x,t;\sita)&:=\sum\limits_{j: (\boldsymbol x_j,t_j)\in P_I} \loss_{\mathrm{IC},j},\ \loss_{\mathrm{IC},j}:=\dfrac{1}{|P_I|}|\widehat{u}_{\sita}(\boldsymbol x_j,0)-u_0(\boldsymbol x_j)|^2;\\
\loss_{\mathrm{BC}}(\boldsymbol x,t;\sita)&:=\sum\limits_{j: (\boldsymbol x_j,t_j)\in P_B} \loss_{\mathrm{BC},j},\ \loss_{\mathrm{BC},j}=|\widehat u_{\sita, j}-u_{B}(\boldsymbol x_j,t)|^2;
\label{Loss_PINN}
\end{aligned}
\end{equation}
where $\mathcal P(\cdot)$ is defined in Eq. (\ref{pde}), $\widehat u_{\sita, j}:=\widehat u_{\sita}(\boldsymbol x_j, t_j)$ denotes the predicted $u$ at $(\boldsymbol x_j, t_j)$; $P_N \subseteq \mathrm{int}\ \Omega \times (0,T], P_I \subseteq \Omega \times \{0\}$ and $P_B \subseteq \partial\Omega \times(0,T]$ denote the collocation point sets in the spatiotemporal solution domain interior, on the initial/boundary condition surface respectively, with $T$ being the terminal time. The weighting coefficients $w_{\mathrm{GOV}}, w_{\mathrm{IC}}, w_{\mathrm{BC}}$ balance the contributions of the governing equation term and initial/boundary condition term in the loss function. Noting that Dirichlet boundary conditions are imposed in this paper, i.e.
\begin{equation}
    u(\boldsymbol x,t)=u_B(\boldsymbol x,t),\quad \forall(\boldsymbol x,t)\in\partial \Omega\times(0,T].
\end{equation}

The parameters in the neural network are then iteratively updated through optimization. Let $\boldsymbol\theta_j$ denote the parameters after $j$ updates, then the parameter update rule can be universally expressed as
\begin{equation}
\sita_{j+1}=OPT({\sita}_j, \nabla_{\sita}\mathcal  L(x,t;\boldsymbol{\theta_j}), \alpha_j),
\label{opt}
\end{equation}
where $OPT$ represents the optimization strategy and $\alpha_j > 0$ denotes the iteration-dependent learning rate. For instance, in Adam, which is widely adopted in PINN, the update takes the specific form
\begin{gather}
\sita_{j+1}=\sita_j-\frac{\alpha_j\boldsymbol{m}_j}{\sqrt{\boldsymbol{v}_j}+\varepsilon};\notag\\
\boldsymbol{m}_{j}=\frac{\beta_1\boldsymbol{m}_{j-1}+(1-\beta_1)\nabla_{\sita} \loss({\boldsymbol{x};\boldsymbol\theta_j})}{1-\beta_1^j};\\
\boldsymbol{v}_{j}=\frac{\beta_2\boldsymbol{v}_{j-1}+(1-\beta_2)\big(\nabla_{\sita} \loss({\boldsymbol{x};\boldsymbol\theta_j})\big)^2}{1-\beta_2^j},\notag
\end{gather}
where $\beta_1,\beta_2\in(0,1)$ are constant, $\boldsymbol{m}_0=\boldsymbol{v}_0=\bf 0$, and $\varepsilon>0$ is a small quantity approaching zero.

\section{Methodology}
\label{sec:method}

In this section, we present the framework of CLINN, which encodes conservation law features into PINN. Then introduce two key techniques for handling discontinuities: the AN-based shock indicator for discontinuity identification and the improved RAR method for discontinuity treatment.

\begin{figure}[H] 
    \centering
\includegraphics[width=\textwidth]{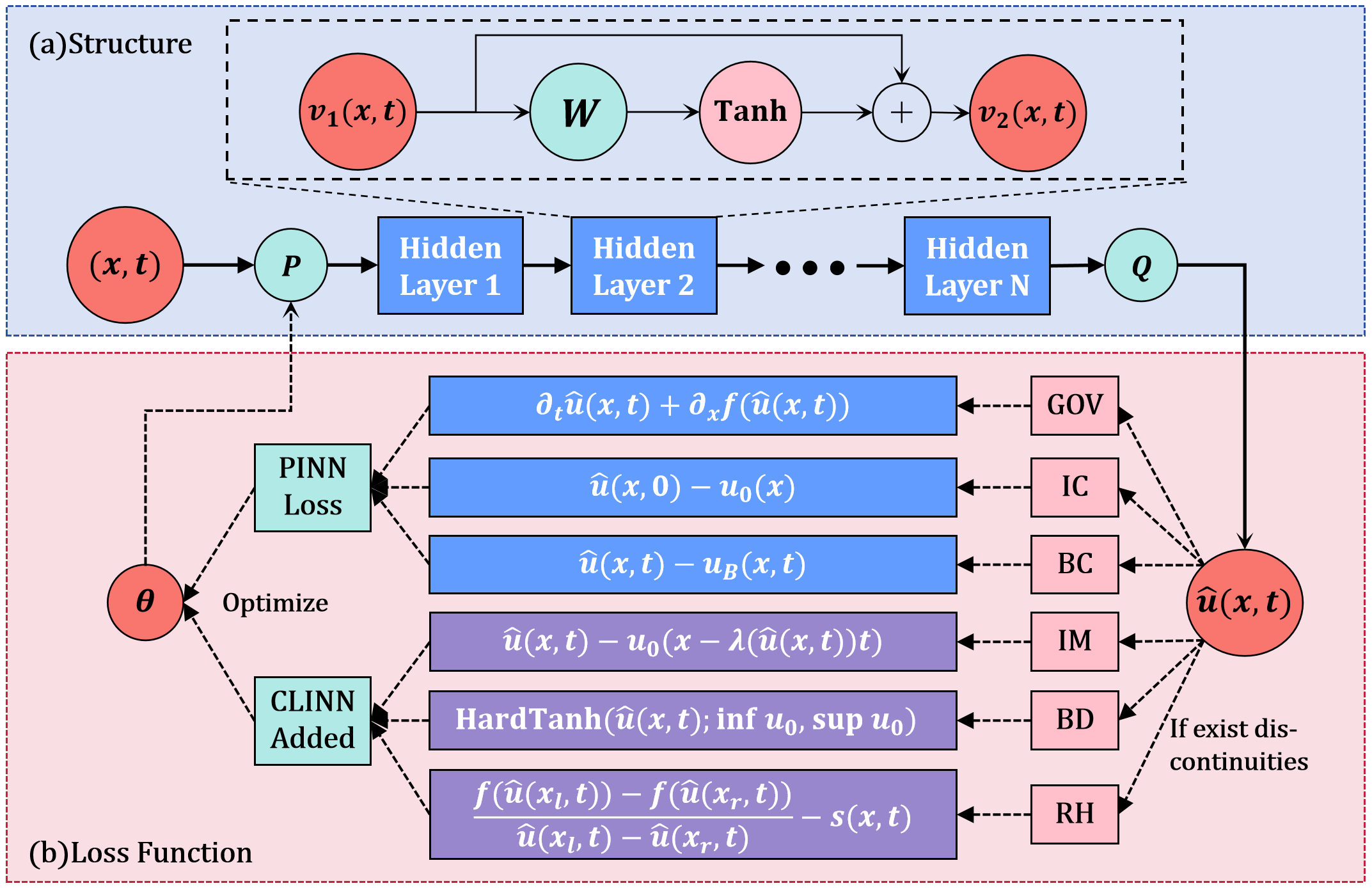}
    \caption{The architecture of CLINN for one-dimensional problems. (a) Structure; (b) Loss function}
    \label{clinn} 
\end{figure}

\subsection{Conservation law informed neural network (CLINN)}

To enhance the convergence speed of PINN while mitigating gradient vanishing/explosion issues, we employ the architecture illustrated in Fig. \ref{clinn}(a). The workflow proceeds as follows. To begin with, coordinate points $(\boldsymbol x,t)$ undergo affine dimension elevation through layer $P: \mathbb{R}^{d+1} \rightarrow \mathbb{R}^n$, producing vector $v_0(\boldsymbol  x,t) = P([\boldsymbol x; t])$ where $n$ denotes the number of neurons in the hidden layer. Subsequent processing through $N$ hidden layers yields $v_N(\boldsymbol x,t) \in \mathbb{R}^n$, with each layer implementing skip connections \cite{he2016deep} formulated as $v_{k+1}(\boldsymbol x,t) = v_k(\boldsymbol x,t) + \mathrm{Tanh}(W_k(v_k(\boldsymbol x,t))$ for $k \in [0, N-1]$, where $W_k: \mathbb{R}^n \rightarrow \mathbb{R}^n$ represents affine transformation. The final output $\widehat{u}(\boldsymbol x,t) = Q(v_N(\boldsymbol x,t))$ is obtained via affine dimension reduction layer $Q: \mathbb{R}^n \rightarrow \mathbb{R}$.

As a general-purpose PDE-solving framework, the loss functions in PINN fail to incorporate distinctive characteristics of conservation laws, which leads to poor performance in handling problems involving discontinuities. To address this limitation, we extend the loss function in IFNN \cite{zhang2022implicit} and PINN-WE \cite{liu2024discontinuity}, resulting in the following enhanced loss construction:
\begin{equation}
\begin{aligned}
\loss(\boldsymbol x,t;\sita)&:=\loss_{\mathrm{PINN}}(\boldsymbol x,t;\sita)+\loss_{\mathrm{CHAR}}(\boldsymbol x,t;\sita);\\
\loss_{\mathrm{CHAR}}(\boldsymbol x,t;\sita)&:=
w_{\mathrm{IM}}\loss_{\mathrm{IM}}(\boldsymbol x,t;\sita)+w_{\mathrm{BD}}\loss_{\mathrm{BD}}(\boldsymbol x,t;\sita)+w_{\mathrm{RH}}\loss_{\mathrm{RH}}(\boldsymbol x,t,\sita);\\
\loss_{\mathrm{IM}}(\boldsymbol x,t,\sita)&:=\sum\limits_{j: (\boldsymbol x_j,t_j)\in P_N}\loss_{\mathrm{IM},j},\ \loss_{\mathrm{IM},j}=\dfrac{1}{|P_N|}\left|\widehat{u}_{\sita, j}-u_0\big(\boldsymbol x_j-\boldsymbol\lambda(\widehat{u}_{\sita, j})t_j\big)\right|^2;\\
\loss_{\mathrm{BD}}(\boldsymbol x,t,\sita)&:=\sum\limits_{j: (\boldsymbol x_j,t_j)\in P_N}\loss_{\mathrm{BD},j}, \ \loss_{\mathrm{BD},j}=\dfrac{1}{|P_N|}\mathrm{HardTanh}^2\Big(\widehat{u}_{\sita, j};\inf\, u_0,\sup\, u_0\Big);\\
\loss_{\mathrm{RH}}(\boldsymbol x,t,\sita)&:=\sum\limits_{j: (x_j,t_j)\in P_D}\loss_{\mathrm{RH},j}, \ \loss_{\mathrm{RH},j}=\left|\frac{\nabla_{\boldsymbol x}\Phi(\boldsymbol x_j, t_j)\cdot [\boldsymbol f(\widehat u_{\sita, j})]}{|\nabla_{\boldsymbol x}\Phi(\boldsymbol x_j, t_j)|[\widehat u_{\sita, j}]}-s(\boldsymbol x_j,t_j)\right|,
\label{Loss_Proposed}
\end{aligned}
\end{equation}
where the function $\mathrm{HardTanh}(\cdot)$ and the point set $P_N$ maintains the same definition as Eq. (\ref{act_func}) and Eq. (\ref{Loss_PINN}) respectively, while ``$\inf u_0$" and ``$\sup u_0$" denote the infimum and supremum of the initial value function $u_0(\boldsymbol x)$ over the domain $\Omega$, respectively. $P_D\subset \Omega\times[0,T]$ represents the set of discontinuity points, and the discontinuity surface $\Phi(\boldsymbol x, t)$ can be obtained through interpolation; $s(\boldsymbol x,t)$ is defined in Proposition \ref{prop}. For one-dimensional problem, the term $\loss_{\mathrm{RH}, j}$ can be simplified as
\begin{equation}
\loss_{\mathrm{RH},j}=\left|\frac{f(\widehat{u}_{j,L})-f(\widehat{u}_{j,R})}{\widehat{u}_{j,L}-\widehat{u}_{j,R}}-s(x_j,t_j)\right|, 
\end{equation}
where $\widehat{u}_{j,L}:=\widehat{u}_{\sita}(x_j-h,t_j), \widehat{u}_{j,R}:=\widehat{u}_{\sita}(x_j+h,t_j)$ denote the predicted values of the neural network at the points on either side of the discontinuity, and $h$ is a small increment. The terms $w_{\mathrm{IM}}, w_{\mathrm{BD}}, w_{\mathrm{RH}}$ represent the weighting coefficients for the implicit solution term, boundedness constraint term, and discontinuity jump condition term, respectively.

In summary, the loss function presented in this paper includes three terms related to PDEs (governing equation, initial/boundary conditions) and three characteristic terms related to conservation laws (implicit solution, boundedness constraint, discontinuity jump condition), as shown in Fig. \ref{clinn}(b).

\subsection{Discontinuity indicator based on artificial neural (AN)}

As demonstrated in Eq. (\ref{Loss_Proposed}), our proposed methodology integrates discontinuity-related features directly into the loss function formulation. This approach naturally raises a fundamental question regarding the discontinuity  identification. Building upon the work of Feng et al. \cite{feng2020characteristic, feng2021characteristic}, who developed an artificial neuron (AN)-based shock indicator that encodes the mathematical properties of conservation laws within a computationally efficient neural network architecture, we can easily capture the discontinuity. Noting that the AN-based indicator overcomes the limitations of traditional shock detection techniques, such as the TVB troubled-cell indicator \cite{cockburn1989tvb}, which incorrectly identifies extrema and large gradients as troubled cells.

For one-dimensional problem, let the solution interval for the scalar conservation law be denoted as $[p,q]$, which is partitioned into $\bigcup_{j=1}^N I_j$, where $I_j=[x_{j-\frac12}, x_{j+\frac12}]$ and satisfies $p=x_\frac12<x_\frac32<\cdots<x_{N+\frac12}=q$. Additionally, $x_j=\frac12(x_{j-\frac12}+x_{j+\frac12})$ is defined as the center point of $I_j$, and $h_j=x_{j+\frac12}-x_{j-\frac12}$ represents the size of the interval $I_j$.

\begin{figure}
    \centering
    \includegraphics[width=\linewidth]{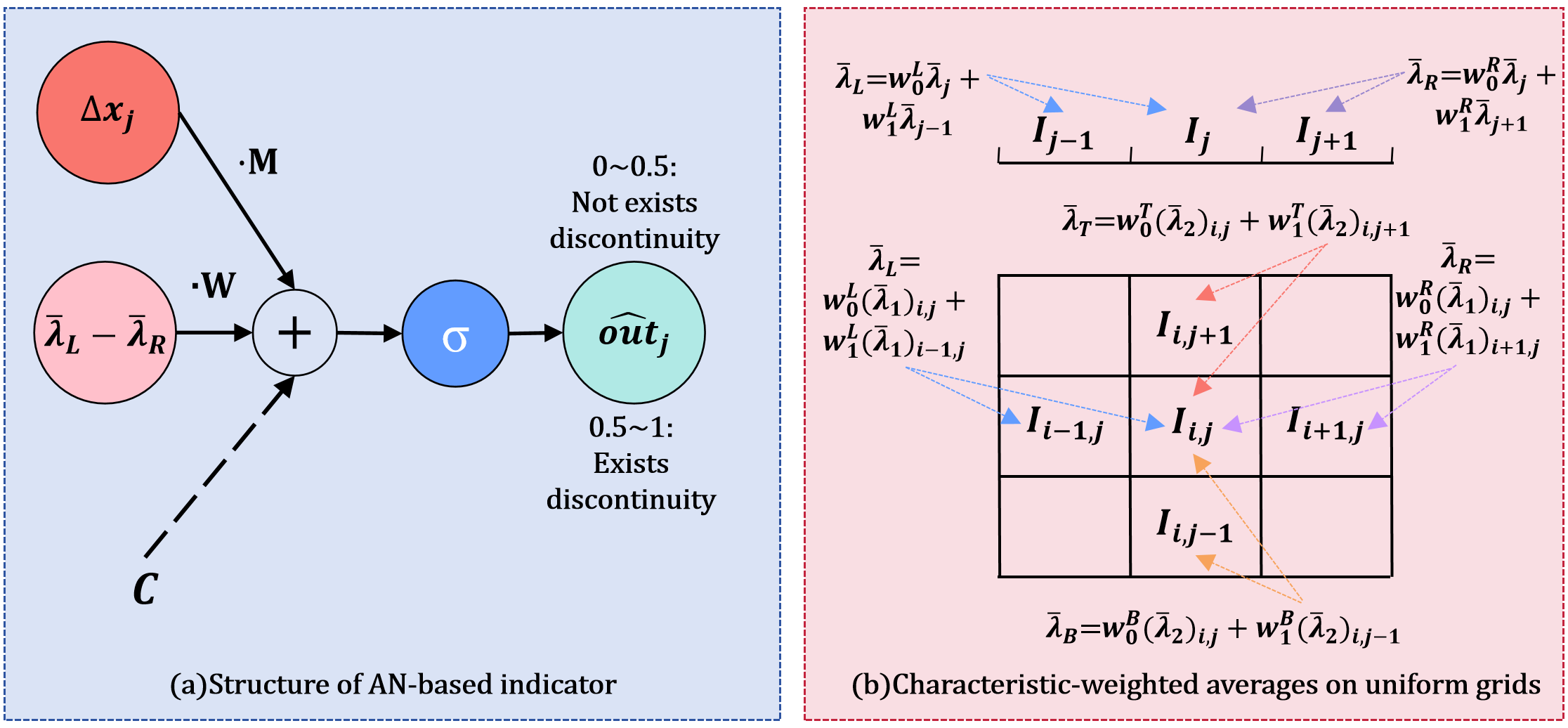}
    \caption{The principle of the AN-based indicator. (a) The structure of the indicator; (b)The characteristic-weighted averages on uniform grids.}
    \label{an_image}
\end{figure}

Feng et al. posit that the existence of a discontinuity within cell $I_j$ at time $t$ is determined by the mesh size $\Delta x_j:=\max\{h_{j-1},h_j,h_{j+1}\}$ and characteristic integral averages $\overline{\lambda}_{j-1}$, $\overline{\lambda}_{j}$, $\overline{\lambda}_{j+1}$, where $\overline\lambda_{j}$ is defined as
\begin{equation}
\overline\lambda_{j}=\frac{1}{h_j}\int_{x_{j-\frac12}}^{x_{j+\frac12}}\lambda(u(\cdot,t))\de x.
\end{equation}
Then construct a FFN (shown in Fig. \ref{an_image} (a)) with only two neurons to detect discontinuities:
\begin{equation}
\widehat{out}_j(\overline{\lambda}_{L}, \overline{\lambda}_{R}, \Delta x_j):=\mathrm{Sigmoid}\left(W(\overline{\lambda}_{L}-\overline{\lambda}_{R})+M\Delta x_j+C\right),
\label{AN}
\end{equation}
where $\overline{\lambda}_{L}, \overline{\lambda}_{R}$ are defined as
\begin{equation}
\begin{cases}
\overline{\lambda}_{L}:=\dfrac{h_{j-1}}{h_{j-1}+h_j}\overline{\lambda}_{j-1}+\dfrac{h_{j}}{h_{j-1}+h_j}\overline{\lambda}_{j}, \\ \overline{\lambda}_{R}:=\dfrac{h_{j+1}}{h_{j+1}+h_j}\overline{\lambda}_{j+1}+\dfrac{h_{j}}{h_{j+1}+h_j}\overline{\lambda}_{j}.

\end{cases}
\label{ax_lambda}
\end{equation}

Based on the results of Feng et al., the parameters of the neural network should be set as $$W=10.0, M=-12.0, C=-1.0, $$and one can determine whether discontinuities exist by the following facts:

\begin{itemize}
    \item 
    If the solution is smooth in $I_j$, or $I_j$ is near the smooth extrema or admissible large gradients, then $\widehat{out}_j<0.5$;
    \item 
    If there is a discontinuity within $I_j$, then $\widehat{out}_j>0.5$.
\end{itemize}

In this paper, employing the midpoint formula instead of the characteristic integral average can achieve good effects in discontinuity indication in this work: 
\begin{equation}
\frac{1}{h_j}\int_{x_{j-\frac12}}^{x_{j+\frac12}}\lambda(u(\cdot,t))\de x = \lambda(u(x_j, t)) + O(\Delta x_j^2).
\label{medium-formula}
\end{equation}
 
Notably, for double contact discontinuities, the characteristic speed difference between adjacent grid cells, $\overline{\lambda}_{L} - \overline{\lambda}_{R}$, vanishes identically, thereby rendering the indicator ineffective in such cases. This observation demonstrates that while the AN-based indicator successfully detects various discontinuity types—including shock waves and left/right contact discontinuities—it remains incapable of identifying double contact discontinuities.

Since equations in high-dimensional space can still be treated as one-dimensional along any given direction, Feng et al. extend the shock indicator for one-dimensional scalar equations to higher dimensions via dimensional splitting. For two-dimensional uniformed grids ${I_{i,j}=[x_{i-\frac12 }, x_{i+\frac12}]\times [y_{j-\frac12}, y_{j+\frac12}]}$, denote $h_i:=x_{i+\frac12}-x_{i-\frac12}$ and $h_j:=y_{j+\frac12}-y_{j-\frac12}$, then the expression of the AN indicator is given by
\begin{equation}
    \widehat{out}_{i,j}=\max\{\widehat{out}_{i}(\overline\lambda_L, \overline\lambda_R, \Delta x_i), \widehat{out}_{j}(\overline\lambda_B, \overline\lambda_T, \Delta y_j)\}
\end{equation}
where $\Delta x_i:=\max\{h_{i-1}, h_i, h_{i+1}\}$ and $\Delta y_j=\max\{h_{j-1}, h_j, h_{j+1}\}$ are defined analogously to the one-dimensional case, and $\overline\lambda_L$, $\overline\lambda_R$, $\overline\lambda_B$, $\overline\lambda_T$ are defined as
\begin{equation}
\begin{cases}
    \overline \lambda_X:=\dfrac{h_{i+k}}{h_{i+k}+h_i}(\overline \lambda_1)_{i+k,j}+\dfrac{h_{i}}{h_{i+k}+h_i}(\overline \lambda_1)_{i,j},\ & k=\pm1\ \mathrm{for}\ X=R/L;\\ 
    \overline \lambda_Y:=\dfrac{h_{j+k}}{h_{j+k}+h_j}(\overline\lambda_2)_{i,j+k}+\dfrac{h_{j}}{h_{j+k}+h_j}(\overline\lambda_2)_{i,j},\ & k=\pm1\ \mathrm{for}\ Y=T/B;
    
\end{cases}
\label{an2d}
\end{equation}
which is shown in Fig. \ref{an_image}(b). Denote $\boldsymbol\lambda = (\lambda_1, \lambda_2)$ and $x_i:=\frac12({x_{i-\frac12}+x_{i+\frac12}}), y_j:=\frac12({y_{j-\frac12}+y_{j+\frac12}})$, then $(\overline{\lambda}_1)_{i,j}$ and $(\overline{\lambda}_2)_{i,j}$ in Eq. (\ref{an2d}) are defined as
\begin{equation}
(\overline{\lambda}_1)_{i,j}=\int_{x_{i-\frac12}}^{x_{i+\frac12}}\lambda_1(u(\cdot, y_j,t))\de x,\quad (\overline{\lambda}_2)_{i,j}=\int_{y_{j-\frac12}}^{y_{j+\frac12}}\lambda_2(u(x_i, \cdot,t))\de y.
\end{equation}
For practical applications in this paper, such as the discontinuity indicator for the 2D Burgers equation in Section \ref{2D}, we similarly employ the midpoint formula approximation as in Eq. (\ref{medium-formula}).

The extension of the AN indicator to higher dimensions follows a methodology analogous to that employed for the two-dimensional case, demonstrating the strong generalization capability of the indicator.

\subsection{Residual-based adaptive refinement (RAR)}
Approximating discontinuous functions with neural networks has long been one of the challenging research problems. This paper combines approaches from DeepXDE \cite{lu2021deepxde} and GA-PINN \cite{ferrer2024gradient} to propose an improved RAR method suitable for our study:

\begin{itemize}
\item After preliminary training of CLINN, we identify points with the largest errors in the governing equation term $\loss_{\mathrm{GOV}}$ in Eq. (\ref{Loss_PINN}) and implicit solution term $\loss_{\mathrm{IM}}$ in Eq. (\ref{Loss_Proposed}) within non-discontinuous regions, then increase their error weights to enhance the neural network's approximation near discontinuities.
\item After identifying discontinuity locations, we eliminate the error weights of the governing equation term and implicit solution term at these points, computing only the error of the discontinuity jump condition term, aiming to make the neural network focus on optimizing solutions in smooth regions.

\end{itemize}

The pseudocode for CLINN equipped with RAR is presented in Algorithm \ref{algo}.

\begin{algorithm}[H]
\caption{CLINN equipped with RAR}
\label{algo}
\KwIn{Point set $P:=\{(\boldsymbol x_j, t_j)\}_{j=1}^J$ sampled in the solution domain; Learning rate $\alpha$; RAR weights for governing equation term and implicit solution term $w_{EQ}, w_{IF}$}
\KwOut{Neural network's predicted solution $\{\widehat{u}_{\sita}(\boldsymbol x_j,t_j)\}_{j=1}^J$}

Initialize network parameters $\sita$ and RAR weights $w_j \gets 1$ for each sample point $(\boldsymbol x_j, t_j)$, $\forall 1 \leq j \leq J$.

\For {$n=0,1,\cdots, N_\mathrm{RAR}$}{
  \For {$ep=1,2,\cdots, N_{\mathrm{epoch},n}$}{
    Compute the network's predicted values $\widehat u_j\gets\widehat{u}_{\sita}(\boldsymbol x_j,t_j)$ pointwise;
    
    Compute the PINN loss by Eq. (\ref{Loss_PINN}): $\loss\gets\displaystyle\sum_{j:(\boldsymbol x_j,t_j)\in P_N\backslash P_D}w_jw_{\mathrm{GOV}}\loss_{\mathrm{GOV},j}+\displaystyle\sum_{j:(\boldsymbol x_j,t_j)\in P_I}w_jw_{\mathrm{IC}}\loss_{\mathrm{IC},j}+\displaystyle\sum_{j:(\boldsymbol x_j,t_j)\in P_B}w_jw_{\mathrm{BC}}\loss_{\mathrm{BC},j};$
    
    Compute the CLINN loss by Eq. (\ref{Loss_Proposed}):
    $\loss\gets\loss+\displaystyle\sum_{j:(\boldsymbol x_j,t_j)\in P_N\backslash P_D}w_jw_{\mathrm{IM}}\loss_{\mathrm{IM},j}+\displaystyle\sum_{j:(\boldsymbol x_j,t_j)\in P_N\backslash P_D}w_jw_{\mathrm{BD}}\loss_{\mathrm{BD},j}+\displaystyle\sum_{j:(\boldsymbol x_j,t_j)\in P_D}w_jw_{\mathrm{RH}}\loss_{\mathrm{RH},j};$
    
    Optimize the parameters $\sita\gets OPT(\sita, \nabla_{\sita}\loss, \alpha)$ by Eq. (\ref{opt});
    
    \If {$ep=N_{\mathrm{epoch}}$}
    {$w_j \gets 1 \ (1 \leq j \leq J)$, and select $N_{pt}$ points with the largest errors in the governing equation term and implicit solution term from $P_D$ respectively, then set $w_j \gets w_j + w_{EQ}$ and $w_j \gets w_j + w_{IF}$.}
  }
}
\end{algorithm}

\section{Numerical Experiment}
\label{sec:experiment}
In this section, we conduct numerical studies on a series of problems to demonstrate the effectiveness and flexibility of CLINN. 

\subsection{Experimental Setup}We choose both 1D and 2D equations for our experiment, and detailed descriptions are provided in Subsection \ref{1D} and \ref{2D}. The model is constructed with the PyTorch framework and run on a NVIDIA Tesla T4 card.

\begin{itemize}
    \item \textbf{Data Preparations: }We generate points $\{t_j, \boldsymbol x_j\}_{j=1}^{N_t\times N_x^d}$ for the input of CLINN through a uniform partition of the computational domain. In the training stage, for 1D problems, $N_t=64, N_x=512$; and for 2D problems, $N_t=32, N_x=128$. In the evaluation stage, $N_t=200, N_x=800$.

    \item \textbf{CLINN settings: }We implement CLINN with a 5-hidden-layer architecture with 100 neurons per layer, optimized via the Adam algorithm (learning rate $\alpha=1e-4$). The RAR scheme (see Algorithm \ref{algo}) is applied once with $(w_{EQ}, w_{IF})=(33,16)$ and $N_{pt}=500$, and a two-phase epoch configuration is constructed ($N_\mathrm{epoch,0}=5000$ for initial training, and $N_\mathrm{epoch,1}=5000$ for refinement).  The loss function in Eq. (\ref{Loss_Proposed}) incorporates weighted components with $w_\mathrm{GOV}=1, w_\mathrm{IC}=1000, w_{\mathrm{BC}}=10, w_\mathrm{IM}=w_\mathrm{BD}=10000, w_\mathrm{RH}=100.$ 

    \item \textbf{Comparisons: } To evaluate the advancement of CLINN over the state-of-the-art deep learning methods in solving scalar conservation laws, we compare it to PINN \cite{raissi2019physics}, IFNN \cite{zhang2022implicit} and PINN-WE \cite{liu2024discontinuity}, with loss function components detailed in Table \ref{lossfunc_compare}; and ``CLINN w/o RAR" is the ablation model for evaluating the role of RAR, with $w_{EQ}$ and $w_{IF}$ set to be zero. We save each model that achieves the smallest MSE between the predicted and the exact solutions after 1000 epochs. It should be emphasized that this metric is not used to guide the model training process, as explicitly indicated in Eq. (\ref{Loss_Proposed}).
    
    \item \textbf{Indicators: } For quantitative comparison of methods in Table \ref{lossfunc_compare}, we evaluate the minimum mean squared error (MSE) across the entire solution domain (denoted as MSE\_All), and time-specific MSE values at t = T/8, 3T/8, 5T/8, and 7T/8 (labeled as MSE\_T1, MSE\_T2, MSE\_T3, and MSE\_T4, respectively), where $T$ means the  termination time of simulation. To further compare the performance between CLINN and the baseline PINN, we introduce a statistical metric called the ``improvement ratio", defined as
\begin{equation}
    \bigg(1-\dfrac{\mathrm{MSE\_All(CLINN)}}{\mathrm{MSE\_All(PINN)}}\bigg)\times100\%.
\label{improve}\end{equation}

\end{itemize}

\begin{table}[H]
    \caption{Loss function components in CLINN, IFNN, PINN-WE, PINN}
    \centering
    \begin{tabular}{c|cccccc}
    \hline
         & $\loss_{\mathrm{GOV}} $
         & $\loss_{\mathrm{IC}} $
         & $\loss_{\mathrm{BC}} $
         & $\loss_{\mathrm{IM}}$
         & $\loss_{\mathrm{BD}}$
         & $\loss_{\mathrm{RH}}$
    \\\hline
    CLINN     & \CheckmarkBold 
              & \CheckmarkBold
              & \CheckmarkBold
              & \CheckmarkBold
              & \CheckmarkBold
              & \CheckmarkBold
    \\
    IFNN      & \CheckmarkBold
              & \CheckmarkBold
              & \CheckmarkBold
              & \CheckmarkBold
              & \XSolidBrush
              & \XSolidBrush
    \\
    PINN-WE   & \CheckmarkBold
              & \CheckmarkBold
              & \CheckmarkBold
              & \XSolidBrush
              & \XSolidBrush
              & \CheckmarkBold
    \\
    PINN      & \CheckmarkBold
              & \CheckmarkBold
              & \CheckmarkBold
              & \XSolidBrush
              & \XSolidBrush
              & \XSolidBrush
    \\\hline
    \end{tabular}

    \label{lossfunc_compare}
\end{table}

\subsection{One-Dimensional Problems}
\label{1D}We select distinct numerical fluxes, and conduct experiments for each flux under different initial conditions to comprehensively evaluate CLINN's performance.

\subsubsection{Inviscid Burgers equation}
\label{4.1}
We first consider the inviscid Burgers equation \cite{burgers1995mathematical}, with the numerical flux given by
\begin{equation}
    f(u)=\frac{u^2}{2}.
\label{Burgers_flux}
\end{equation}
This problem frequently appears in the literature as a benchmark for evaluating the performance of numerical methods. Experiments are conducted under a periodic initial condition and a piecewise smooth initial condition:
\begin{gather}
    \mathrm{(1A)}\ u_0(x)=\sin(\pi x)+0.5;\\\mathrm{(1B)}\ u_0(x)=3x\cdot\mathbb{I}_{[-1,3]}(x).
\end{gather}
For the test case (1A), the solution is implicitly determined by 
\begin{equation}
     u_{\mathrm{1A}}(x,t)=\sin\left(\pi(x-u_{\mathrm{1A}}(x,t)t)\right)+0.5,
\end{equation}
and requires an iterative method for numerical computation. The trajectory of the discontinuity is given by 
\begin{equation}
    \Phi(x,t)=x-\frac t2-(2k+1)=0, \quad t>\frac1\pi,\ x\in\mathbb Z.
\end{equation}

\begin{figure}[H]
    \centering
    \includegraphics[width=\linewidth]{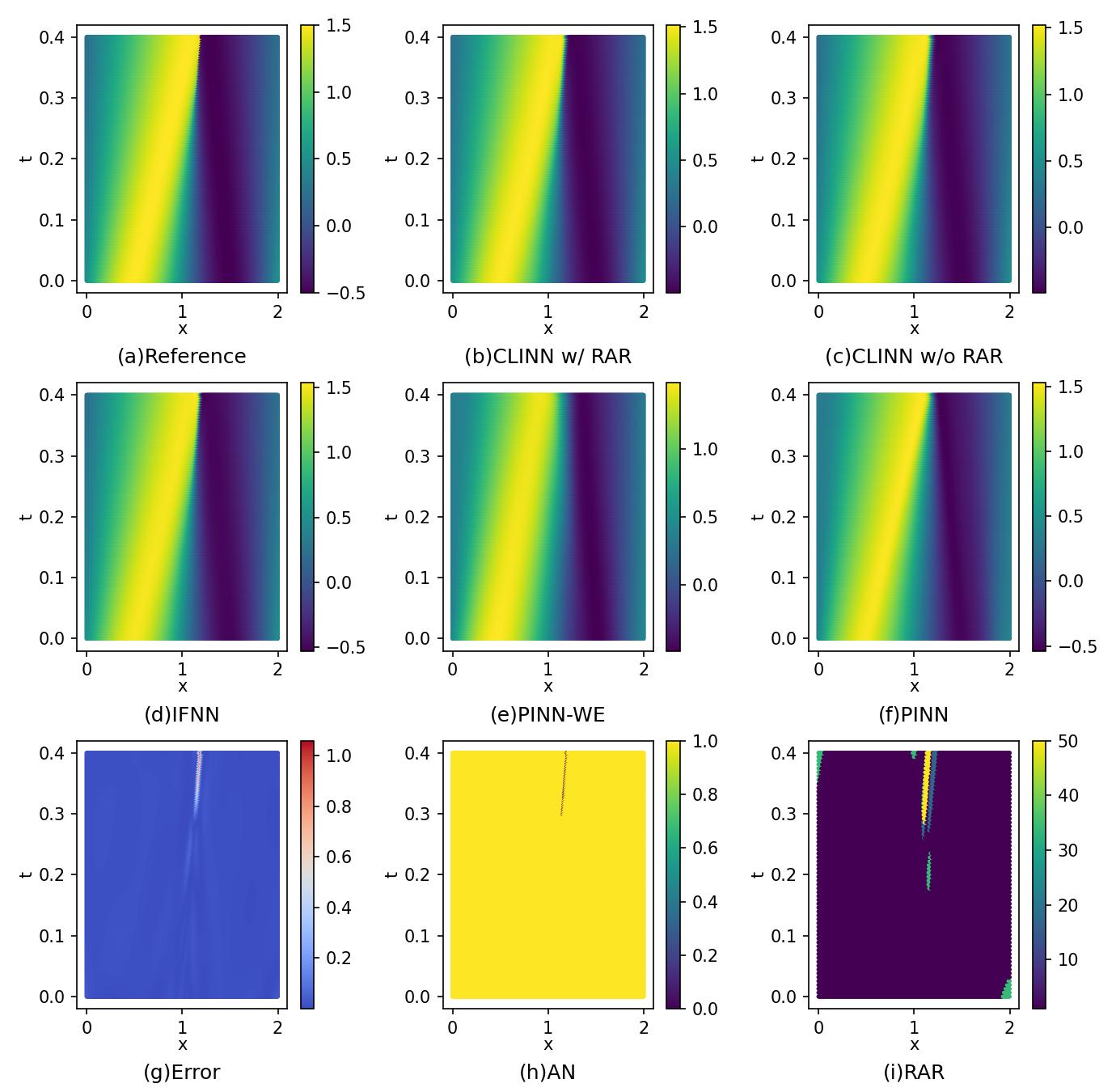}
    \includegraphics[width=\linewidth]{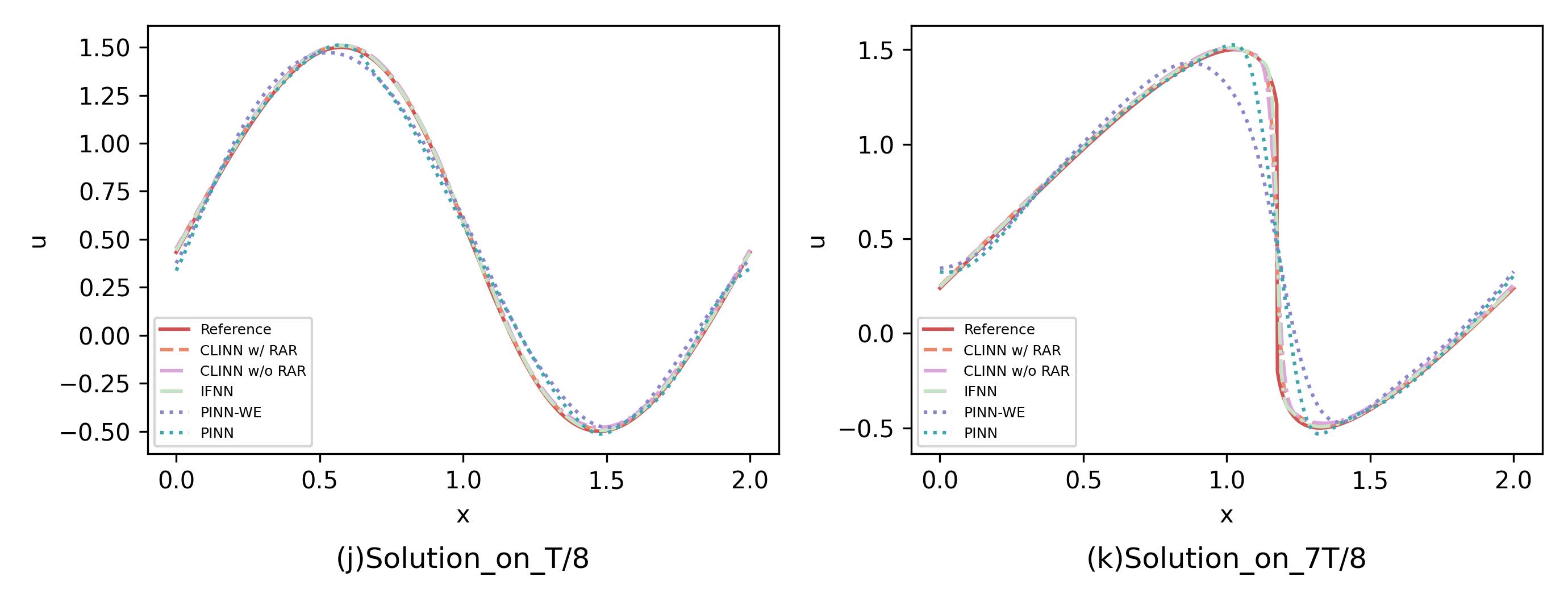}
    \caption{Performance of all methods in test case (1A). (a) Reference solution $(u_{\mathrm{true}})$; (b) Result of CLINN with RAR $(u_{\mathrm{CLINN}})$; (c) Result of CLINN without RAR; (d) Result of IFNN; (e) Result of PINN-WE; (f) Result of PINN; (g) Error of CLINN $(|u_{\mathrm{CLINN}}-u_{\mathrm{true}}|)$; (h) Discontinuity indicator result of AN; (i) Weight distribution in RAR; (j) Comparative results of all methods at $t=T/8$; (k)Comparative results of all methods at $t=7T/8$.}
    \label{case_2}
\end{figure}

\begin{figure}[H]
    \centering
    \includegraphics[width=\linewidth]{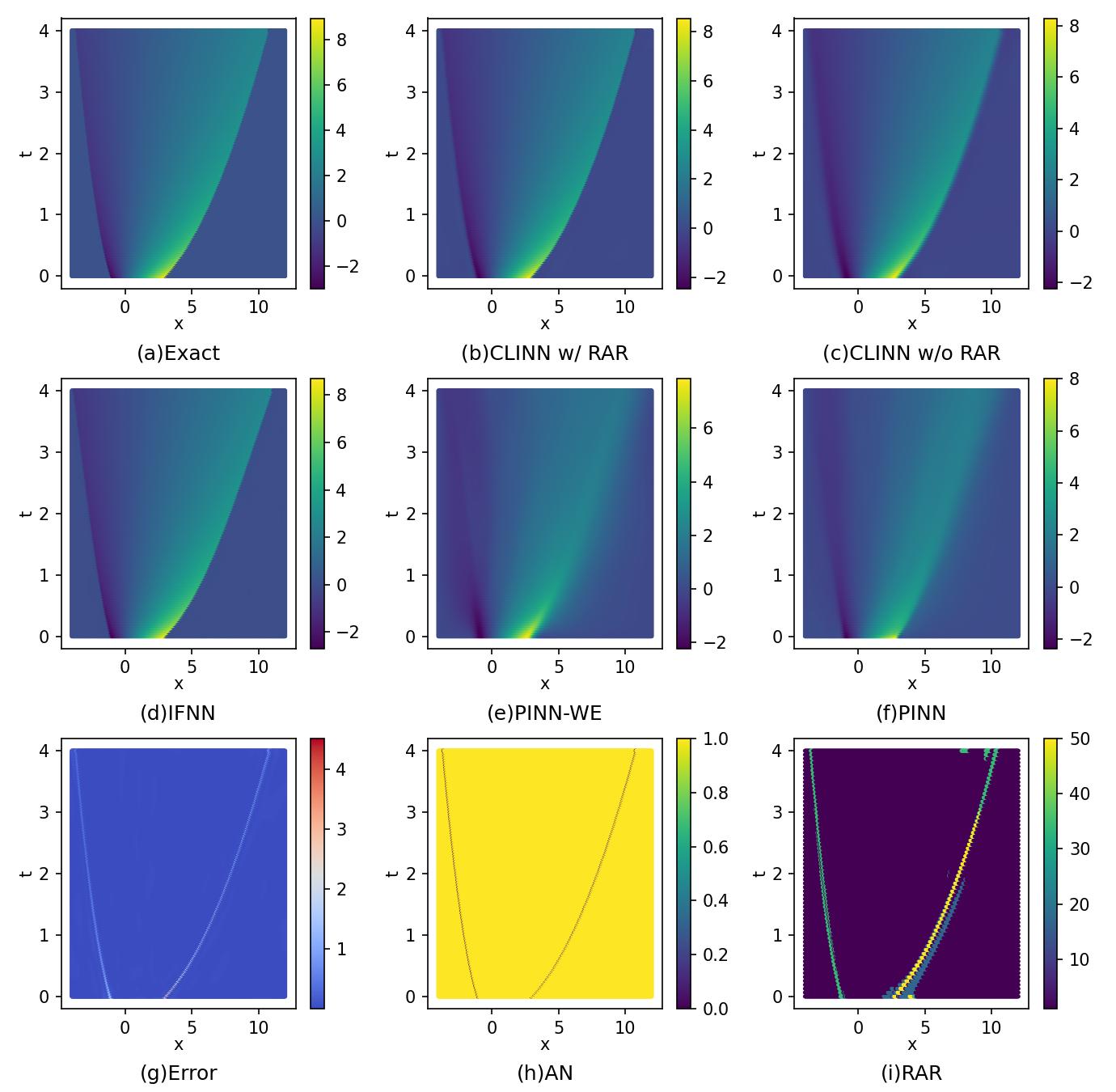}
    \includegraphics[width=\linewidth]{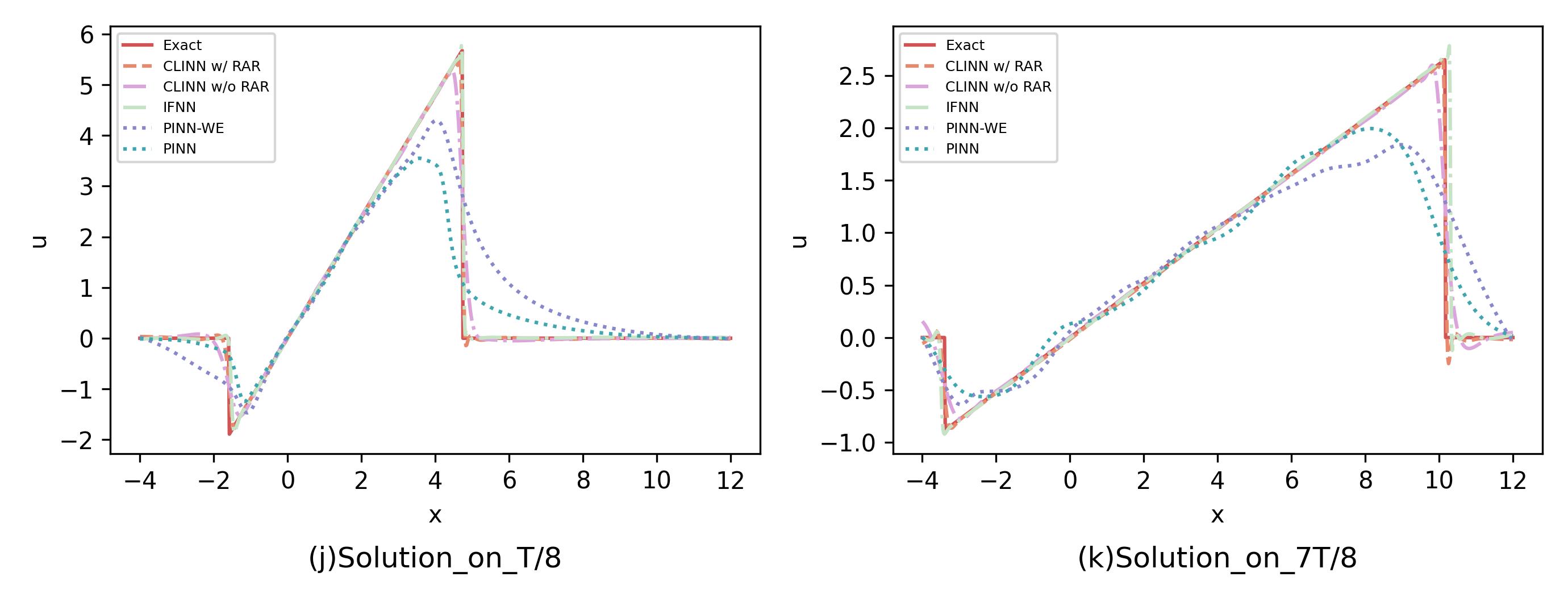}
    \caption{Performance of all methods in test case (1B). (a) Exact solution $(u_{\mathrm{true}})$; (b) Result of CLINN with RAR $(u_{\mathrm{CLINN}})$; (c) Result of CLINN without RAR; (d) Result of IFNN; (e) Result of PINN-WE; (f) Result of PINN; (g) Error of CLINN $(|u_{\mathrm{CLINN}}-u_{\mathrm{true}}|)$; (h) Discontinuity indicator result of AN; (i) Weight distribution in RAR; (j) Comparative results of all methods at $t=T/8$; (k)Comparative results of all methods at $t=7T/8$.}
    \label{case_7}
\end{figure}
For the test case (1B), the exact solution is given by
\begin{equation}
    u_{\mathrm{1B}}(x,t)=
    \begin{cases}\frac{3x}{1+3t}, &-\sqrt{1+3t}\leq x\leq 3\sqrt{1+3t};\\0, &\mathrm{otherwise}.
    \end{cases}
\end{equation}
We solve (1A) on $(x,t)\in[0,2]\times[0,\frac25]$ and (1B) on $[-4,12]\times[0, 4]$. 

The boundary conditions for test case (1A) are specified as
\begin{equation}
    u(0,t)=u(2,t)=g(t), t\in(0,\frac25],
\end{equation}
where $g(t)$ satisfies $g(t)=\sin(-\pi tg(t))+0.5$; and for test case (1B) and all cases in Subsection \ref{4.2} and \ref{4.3}, assuming the computational domain is $(x,t)\in[p,q]\times[0,T]$, then the boundary conditions are given by
\begin{equation}
u(p, t) = u_0(p), \quad u(q, t) = u_0(q), \quad \forall t \in (0, T].
\end{equation}

Sub figures (a)-(g), (j)-(k) in Fig. \ref{case_2} and \ref{case_7} present the performance of all methods on (1A) and (1B). For these test cases, CLINN and IFNN demonstrate marginally superior accuracy compared to PINN-WE and PINN, indicating that incorporating the implicit solution term in the loss function effectively reduces solution errors near discontinuities in this scenario. Sub figure (h) displays the discontinuity detection results of the AN-based indicator, where identified discontinuity points are marked as 0 and smooth regions as 1. The AN-predicted discontinuity locations show excellent agreement with the exact shock positions. Sub figure (i) visualizes the spatiotemporal distribution of RAR weights, revealing that points with elevated weights predominantly cluster near discontinuities. This demonstrates the adaptive refinement's effectiveness in prioritizing critical regions during training.

Table \ref{Compare_MSE_1} presents the error comparison between CLINN and other advanced deep learning methods, demonstrating that for test case (1B), CLINN achieves superior global MSE and time-specific MSE compared to other deep learning methods. However, for test case (1A), CLINN shows slightly inferior performance relative to IFNN, suggesting that the incorporation of jump condition terms may potentially impede the neural network's optimization process and slow down convergence near discontinuities in relatively simpler test cases. The improvement ratios (see Eq. (\ref{improve})) of CLINN over the PINN reach 87.4\% and 96.3\% for these two test cases, demonstrating that the incorporation of both the implicit solution term and jump condition term significantly enhances the neural network's capability to approximate solutions of conservation laws.

\begin{table}[H]
\caption{Comparation among advanced deep learning methods in Burgers equations}\label{Compare_MSE_1}

\centering
\begin{tabular}{cc|cccc|c}
\hline
Case & Model & 
MSE\_T1 &
MSE\_T2 & 
MSE\_T3 &
MSE\_T4 & 
MSE\_All
\\ \hline

\multirow{4}{*}{(1A)} 
& CLINN
& 1.21e-04
& 1.16e-04
& 1.87e-04
& 3.97e-03
& 1.03e-03
\\
& IFNN
& 4.90e-05
& 5.09e-05
& 8.58e-05
& 2.52e-03
& \bf{6.91e-04}
\\
& PINN-WE
& 3.22e-03
& 9.24e-03
& 2.01e-02
& 3.61e-02
& 1.70e-02
\\

& PINN
& 2.52e-03
& 4.85e-03
& 9.33e-03
& 1.85e-02
& 8.16e-03
\\\hline
\multirow{4}{*}{(1B)} 
& CLINN
& 1.87e-02
& 1.07e-02
& 4.97e-03
& 4.23e-03
& \bf{1.15e-02}
\\
& IFNN
& 3.48e-02
& 1.67e-02
& 2.86e-03
& 6.71e-02
& 2.87e-02
\\
& PINN-WE
& 5.14e-01
& 3.74e-01
& 2.33e-01
& 1.43e-01
& 3.04e-01
\\

& PINN
& 5.85e-01
& 3.49e-01
& 2.21e-01
& 1.36e-01
& 3.11e-01
\\\hline

\end{tabular}
\end{table}

\begin{figure}[H]
    \centering
    \includegraphics[width=\linewidth]{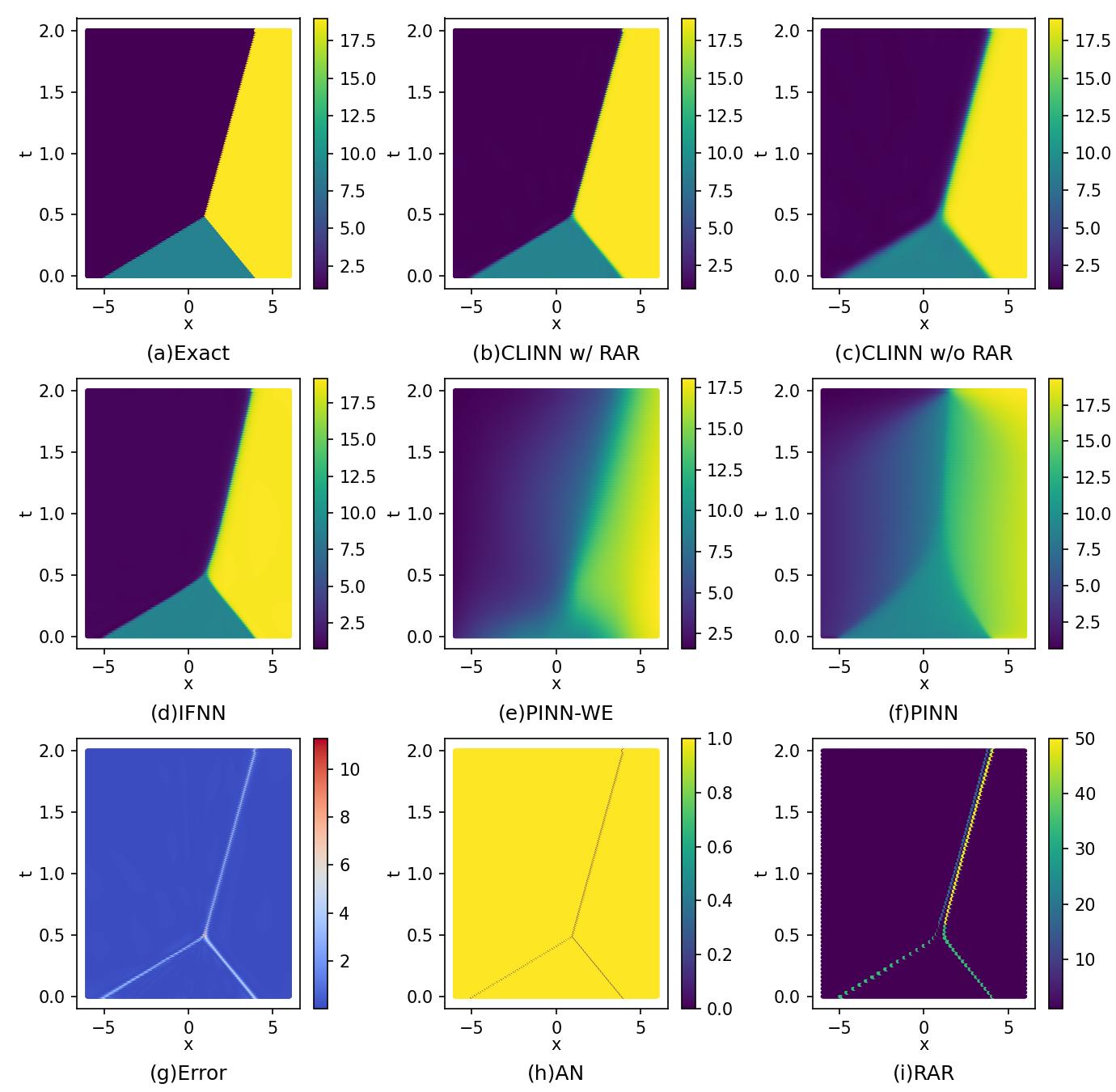}
    \includegraphics[width=\linewidth]{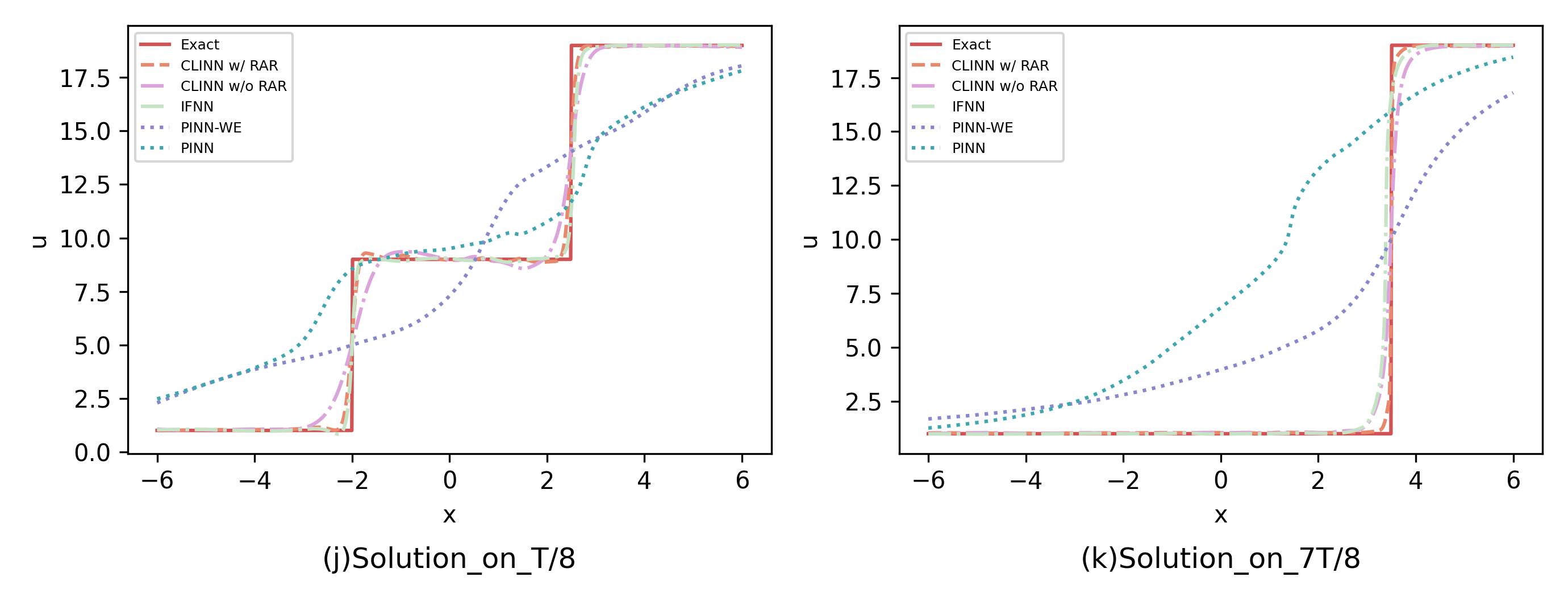}
    \caption{Performance of all methods in test case (2A). (a) Exact solution $(u_{\mathrm{true}})$; (b) Result of CLINN with RAR $(u_{\mathrm{CLINN}})$; (c) Result of CLINN without RAR; (d) Result of IFNN; (e) Result of PINN-WE; (f) Result of PINN; (g) Error of CLINN $(|u_{\mathrm{CLINN}}-u_{\mathrm{true}}|)$; (h) Discontinuity indicator result of AN; (i) Weight distribution in RAR; (j) Comparative results of all methods at $t=T/8$; (k)Comparative results of all methods at $t=7T/8$.}
    \label{case_8}
\end{figure}

\begin{figure}[H]
    \centering
    \includegraphics[width=\linewidth]{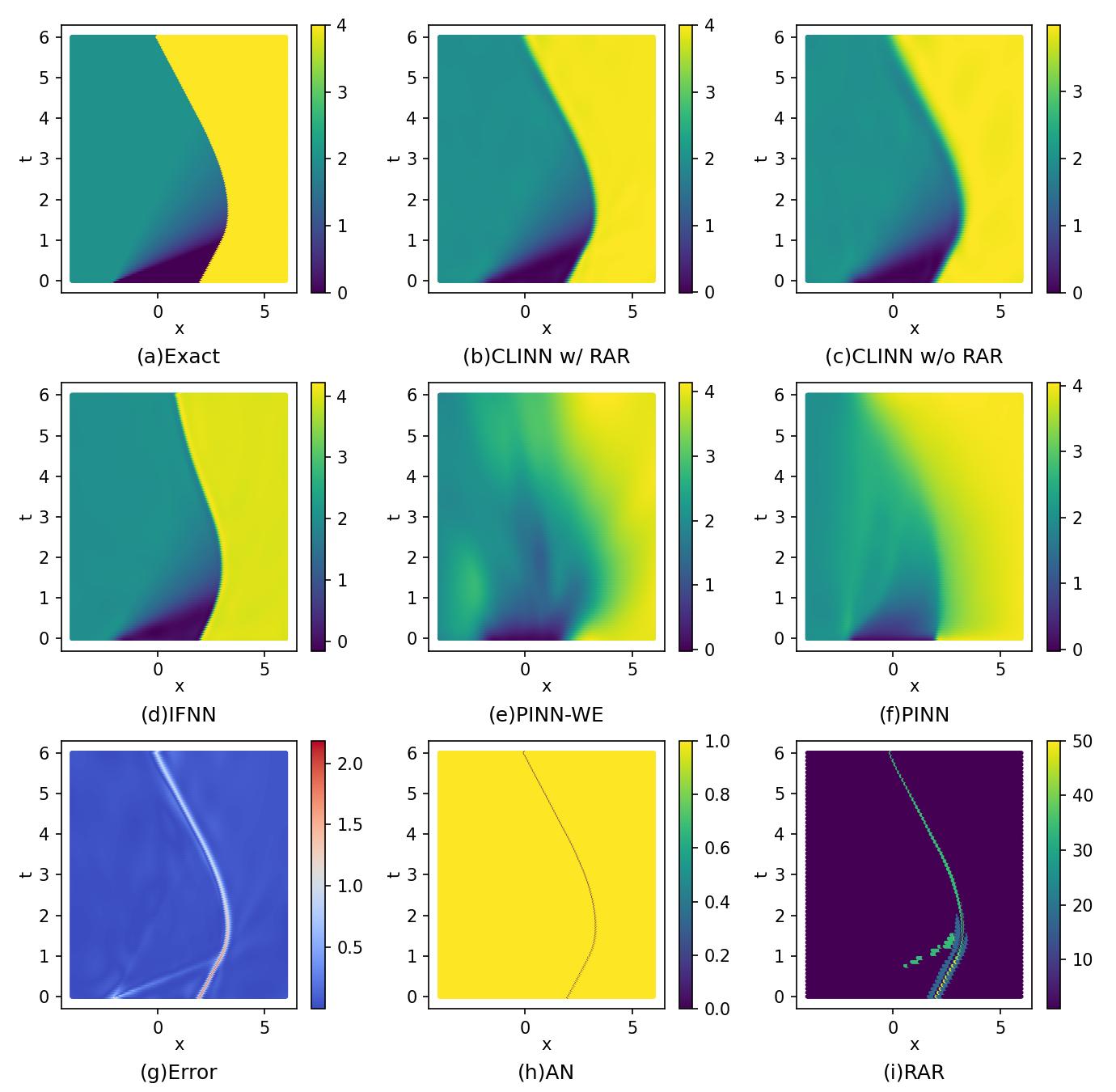}
    \includegraphics[width=\linewidth]{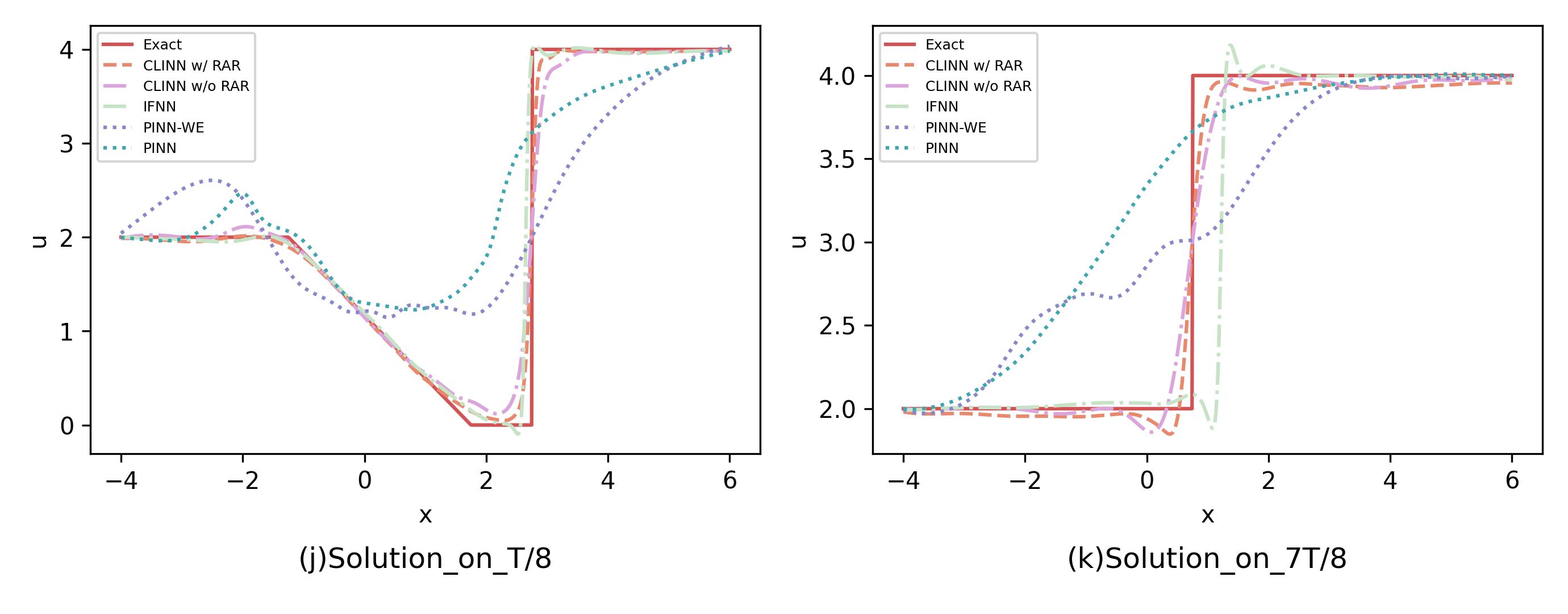}
    \caption{Performance of all methods in test case (2B). (a) Exact solution $(u_{\mathrm{true}})$; (b) Result of CLINN with RAR $(u_{\mathrm{CLINN}})$; (c) Result of CLINN without RAR; (d) Result of IFNN; (e) Result of PINN-WE; (f) Result of PINN; (g) Error of CLINN $(|u_{\mathrm{CLINN}}-u_{\mathrm{true}}|)$; (h) Discontinuity indicator result of AN; (i) Weight distribution in RAR; (j) Comparative results of all methods at $t=T/8$; (k)Comparative results of all methods at $t=7T/8$.}
    \label{case_5}
\end{figure}

\subsubsection{LWR traffic flow model}
\label{4.2}
We now consider Lighthill-Whitham-Richards (LWR) models \cite{lighthill1955kinematic, richards1956shock} in traffic flow theory, with the Greensheld flux \cite{greenshields1935study} given by
\begin{equation}
    f(u)=v_{\max}u\left(1-\frac{u}{u_{\max}}\right),
\end{equation}
where $v, u>0$ represent the velocity and density of the traffic flow. Unlike Eq. (\ref{Burgers_flux}), the flux here is concave. 

Experiments are conducted under two distinct scenarios:
\begin{gather}
    \mathrm{(2A)}\ v_{\max}=u_{\max}=22,\quad u_0(x)=
    \begin{cases}
    1,&x<-5 \\
    9,&-5\leq x<4 \\
    19,&x\geq 4
    \end{cases};\\
    \mathrm{(2B)}\ v_{\max}=u_{\max}=5,\quad u_0(x)=
    \begin{cases}
    2,&x<-2 \\
    0,&-2\leq x<2 \\
    4,&x\geq 2
    \end{cases}.
\end{gather}
The interaction between waves generated from initial discontinuities significantly increases the simulation's complexity. For the test case (2A), where shock waves interact, the solution is given by
\begin{equation}
    u_{\mathrm{2A}}(x,t)=\begin{cases}
    1, & \begin{cases}t<\frac12\ \&\ x<12t-5;\\ t\geq \frac12\ \&\ x<2t;\end{cases}\\ 
    9, & t<\frac12\ \&\ 12t-5\leq x<-6t+4;\\
    19, & \begin{cases}t<\frac12\ \&\ x\geq -6t+4;\\ t\geq \frac12\ \&\ x\geq 2t.\end{cases} \end{cases}
\end{equation}
For the test case (2B), where the centered rarefaction wave interacts with the shock wave, causing the discontinuity's motion to reverse direction, the solution is given by
\begin{equation}
        u_{\mathrm{2B}}(x,t)=\begin{cases}
        2,&
        \begin{cases}
        t<4\ \&\ x<t-2; \\
        t\geq4\ \&\ x>-t+6; 
        \end{cases} \\
        \frac12\left(5-\frac{x+2}{t}\right),&
        \begin{cases}
        t<1\ \&\ t-2\leq x<5t-2; \\
        1\leq t<4\ \&\ t-2\leq x<-2-3t+8\sqrt{t};
        \end{cases} \\
        0,&t<1\ \&\ 5t-2\leq x<t+2; \\
        4, &
        \begin{cases}
        t<1\ \&\ x\geq t+2; \\
        1\leq t<4\ \&\ x\geq -2-3t+8\sqrt{t};
         \\
        t\geq4\ \&\ x\geq -t+6.
        \end{cases}\end{cases}
    \end{equation}
We solve (2A) on $(x,t)\in[-6, 6]\times[0,2]$ and (2B) on $[-4,6]\times[0, 6]$. 

Fig. \ref{case_8} demonstrates the performance of various methods on the test case (2A). The results show that CLINN and IFNN achieve comparable and superior predictive accuracy. Due to the absence of implicit solution constraints, PINN-WE exhibits significant discontinuity smearing - a phenomenon similar to that observed in the test case (1A); while PINN's solutions demonstrate substantial distortion.

Fig. \ref{case_5} compares the performance of different methods on the test case (2B), where the interaction between distinct wave types leads to shifted discontinuity trajectories, resulting in regions where the solution defined by Eq. (\ref{solution}) loses uniqueness and requires supplementary constraints from the jump condition in Eq. (\ref{RH}) to determine the physically admissible solution. A representative example occurs in the region $0\leq 6-t\leq x\leq 2$, where both $u=2$ and $u=4$ satisfy the implicit form relation in Eq. (\ref{solution}), but only $u=4$ is the valid solution through enforcement of the Rankine-Hugoniot condition in Eq. (\ref{RH}). The results of each method further corroborate this observation: without guidance from the jump condition term, IFNN produces discontinuities that deviate from the ground truth when $t > 4$; while PINN-WE, relying solely on the PDE itself and jump conditions but lacking implicit solution constraints, still yields distorted solutions. In contrast, CLINN accurately captures the solution morphology across all spatiotemporal regions through the joint guidance of $\loss_{\mathrm{IM}}$ and $\loss_{\mathrm{RH}}$. Furthermore, as evidenced in Fig. \ref{case_5}(k), the boundedness constraint in Eq. (\ref{bound}) - implemented through the loss term $\mathcal{L}_{\mathrm{BD}}$ - effectively suppresses numerical oscillations on the right-hand side of the discontinuity.

Table \ref{Compare_MSE_2} presents the error of CLINN and other advanced deep learning methods, demonstrating that for the test case (2A) and (2B), CLINN achieves remarkable superiority over all three methods. The improvement ratio of CLINN over PINN reaches 99.2\% and 95.2\% for the two test cases.

\begin{table}[H]
\caption{Comparation among advanced deep learning methods in LWR models}\label{Compare_MSE_2}
\centering
\begin{tabular}{cc|cccc|c}
\hline
Case & Model & 
MSE\_T1 &
MSE\_T2 & 
MSE\_T3 &
MSE\_T4 & 
MSE\_All
\\ \hline

\multirow{4}{*}{(2A)} 
& CLINN
& 3.16e-01
& 1.38e-01
& 1.47e-01
& 1.44e-01
& \bf{2.12e-01}
\\
& IFNN
& 3.98e-01
& 2.26e+00
& 9.27e-01
& 2.12e+00
& 1.40e+00
\\
& PINN-WE
& 9.25e+00
& 1.48e+01
& 1.72e+01
& 1.47e+01
& 1.37e+01
\\

& PINN
& 8.65e+00
& 2.21e+01
& 3.12e+01
& 3.96e+01
& 2.65e+01
\\\hline
\multirow{4}{*}{(2B)} 
& CLINN
& 3.45e-02
& 1.95e-02
& 1.44e-02
& 1.67e-02
& \bf{2.15e-02}
\\
& IFNN
& 1.23e-01
& 1.80e-01
& 3.70e-02
& 1.81e-01
& 1.43e-01
\\
& PINN-WE
& 6.08e-01
& 2.68e-01
& 1.66e-01
& 2.51e-01
& 2.89e-01
\\

& PINN
& 7.68e-01
& 4.78e-01
& 3.37e-01
& 3.33e-01
& 4.49e-01
\\\hline
\end{tabular}
\end{table}

\subsubsection{Non-convex and non-concave flux}
\label{4.3}

All numerical fluxes studied in the previous subsections are either convex or concave, while this subsection investigates non-convex and non-concave fluxes.

We first consider Buckley-Leverett equation \cite{buckley1942mechanism}, with the flux given by
\begin{equation}
    f(u)=\frac{u^2}{u^2+M(1-u)^2}
\end{equation}
where $u \in (0,1)$ represents the water saturation in the oil-water mixture, and $M$ denotes the viscosity ratio between water and oil. In this paper, we set $M = 1$. For the Riemann problem
\begin{equation}
    \mathrm{(3A)}\ u_0(x)=\begin{cases}1, & x\leq1\\0, &x> 1\end{cases},
\end{equation}
the solution develop contact discontinuities, presenting greater computational challenges than problems with convex or concave flux. The solution is given by
\begin{equation}
    u_{\mathrm{3A}}(x,t)=\begin{cases}1, & x\leq 1;\\\frac12(1+\sqrt{1-\frac{t}{x-1}(2\frac{x-1}{t}+1-\sqrt{4\frac{x-1}{t}+1})}), &1< x\leq1+\frac{\sqrt 2+1}{2}t;\\0,&x>1+\frac{\sqrt 2+1}{2}t.\end{cases}
\end{equation}

Another non-convex/non-concave flux we consider is
\begin{equation}
    f(u)=\frac{u^3}{3},
\end{equation}
with the initial value given by
\begin{equation}
    \mathrm{(3B)}\ u_0(x)=
    \begin{cases}
    -2,&x<0 \\
    \frac{3}{2},&0\leq x<1 \\
    1,&x\geq 1
    \end{cases}.
\end{equation}
The current test case contains one more centered rarefaction wave than the previous, with the solution given by
\begin{equation}
    u_{\mathrm{3B}}(x,t)=\begin{cases}
        -2,&x<t;\\
        \sqrt{\frac{x}{t}}, &t\leq x<\frac{9}{4}t;\\\frac32,&\frac94t\leq x\leq \frac94t+1;\\\sqrt{\frac{x-1}{t}},& \frac94t+1< x<4t+1;\\2,&x\geq 4t+1.
    \end{cases}
\end{equation}
We solve (3A) on $(x,t)\in[0, 4]\times[0,2]$ and (3B) on $[-1,3]\times[0, 1]$. 

\begin{figure}[H]
    \centering
    \includegraphics[width=\linewidth]{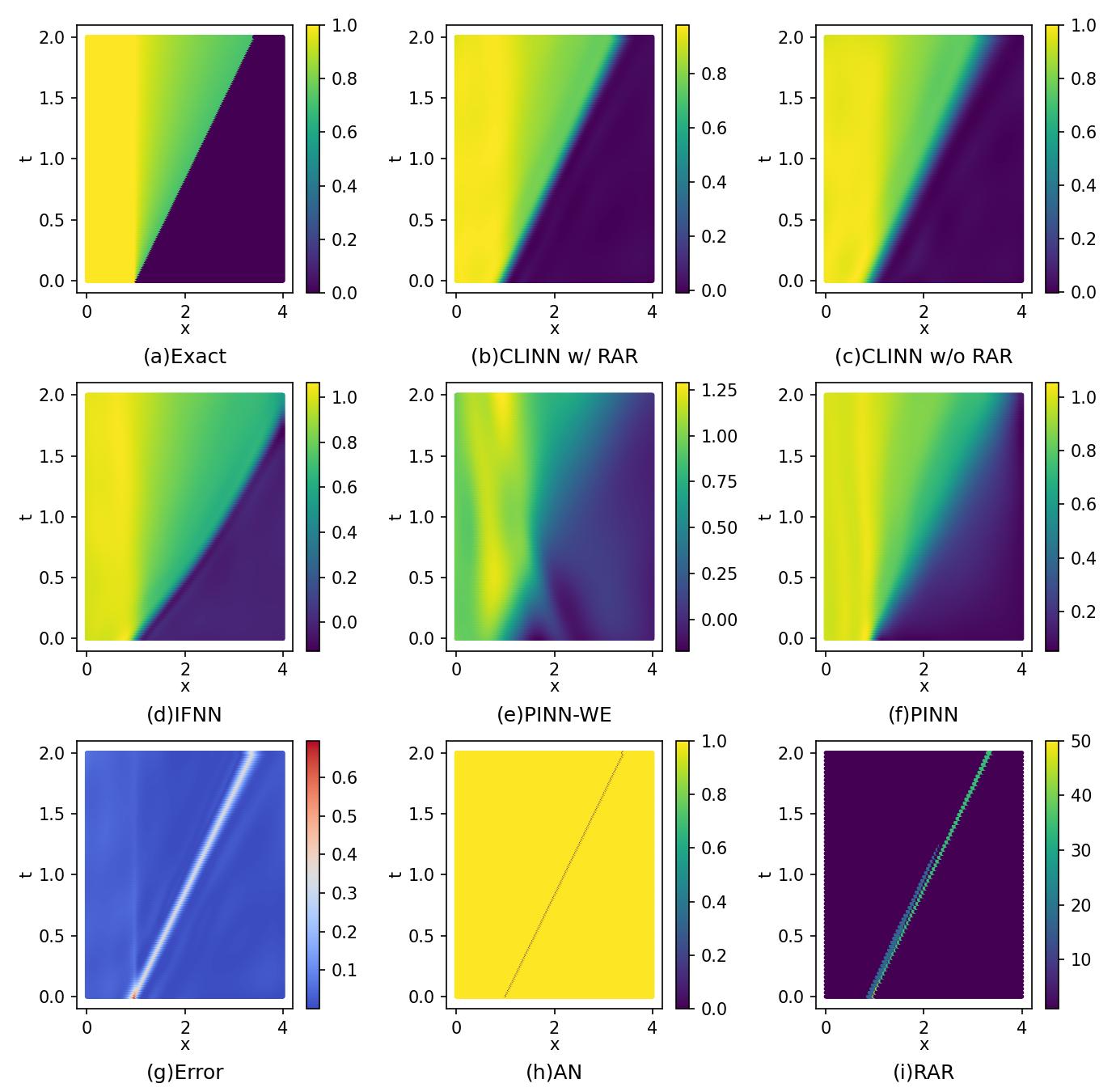}
    \includegraphics[width=\linewidth]{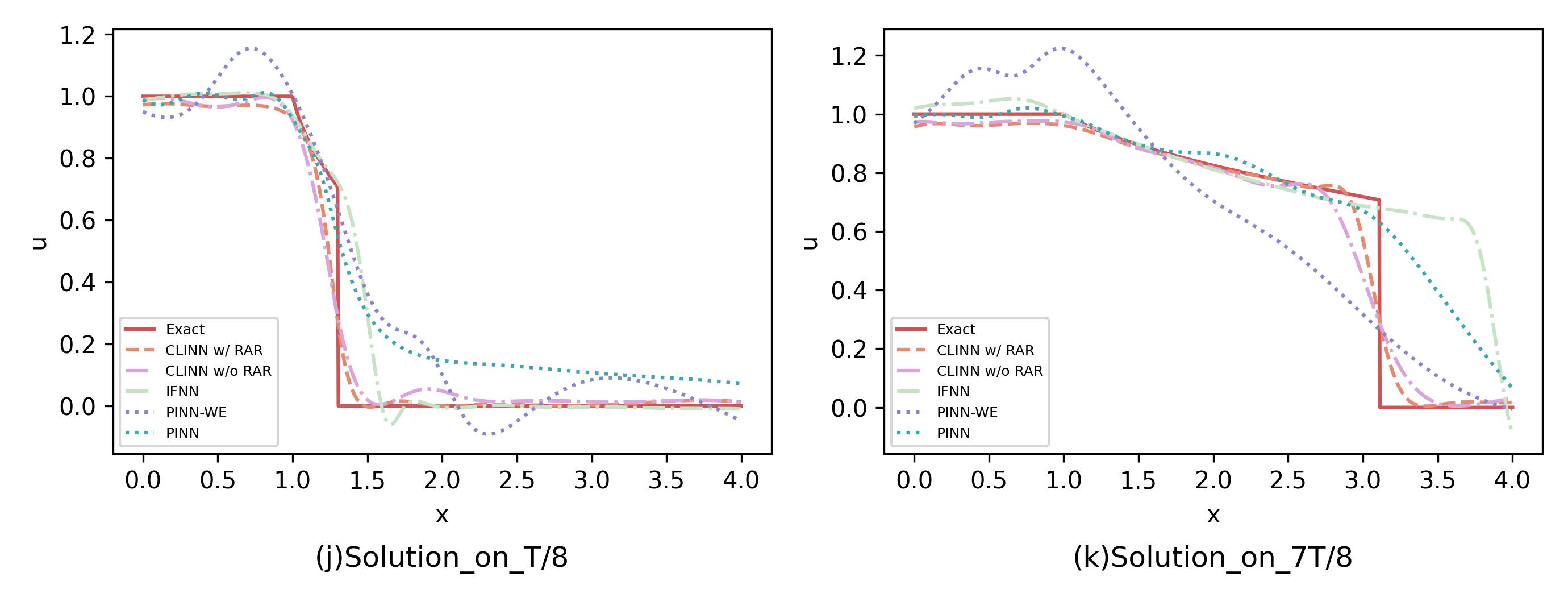}
    \caption{Performance of all methods in test case (3A). (a) Exact solution $(u_{\mathrm{true}})$; (b) Result of CLINN with RAR $(u_{\mathrm{CLINN}})$; (c) Result of CLINN without RAR; (d) Result of IFNN; (e) Result of PINN-WE; (f) Result of PINN; (g) Error of CLINN $(|u_{\mathrm{CLINN}}-u_{\mathrm{true}}|)$; (h) Discontinuity indicator result of AN; (i) Weight distribution in RAR; (j) Comparative results of all methods at $t=T/8$; (k)Comparative results of all methods at $t=7T/8$.}
    \label{case_6}
\end{figure}

Fig. \ref{case_6} demonstrates the performance of different methods on the test case (3A), where CLINN again achieves significantly superior results compared to the other three approaches, both in approximating the centered rarefaction wave region and in predicting the location of the contact discontinuity.

\begin{figure}[H]
    \centering
    \includegraphics[width=\linewidth]{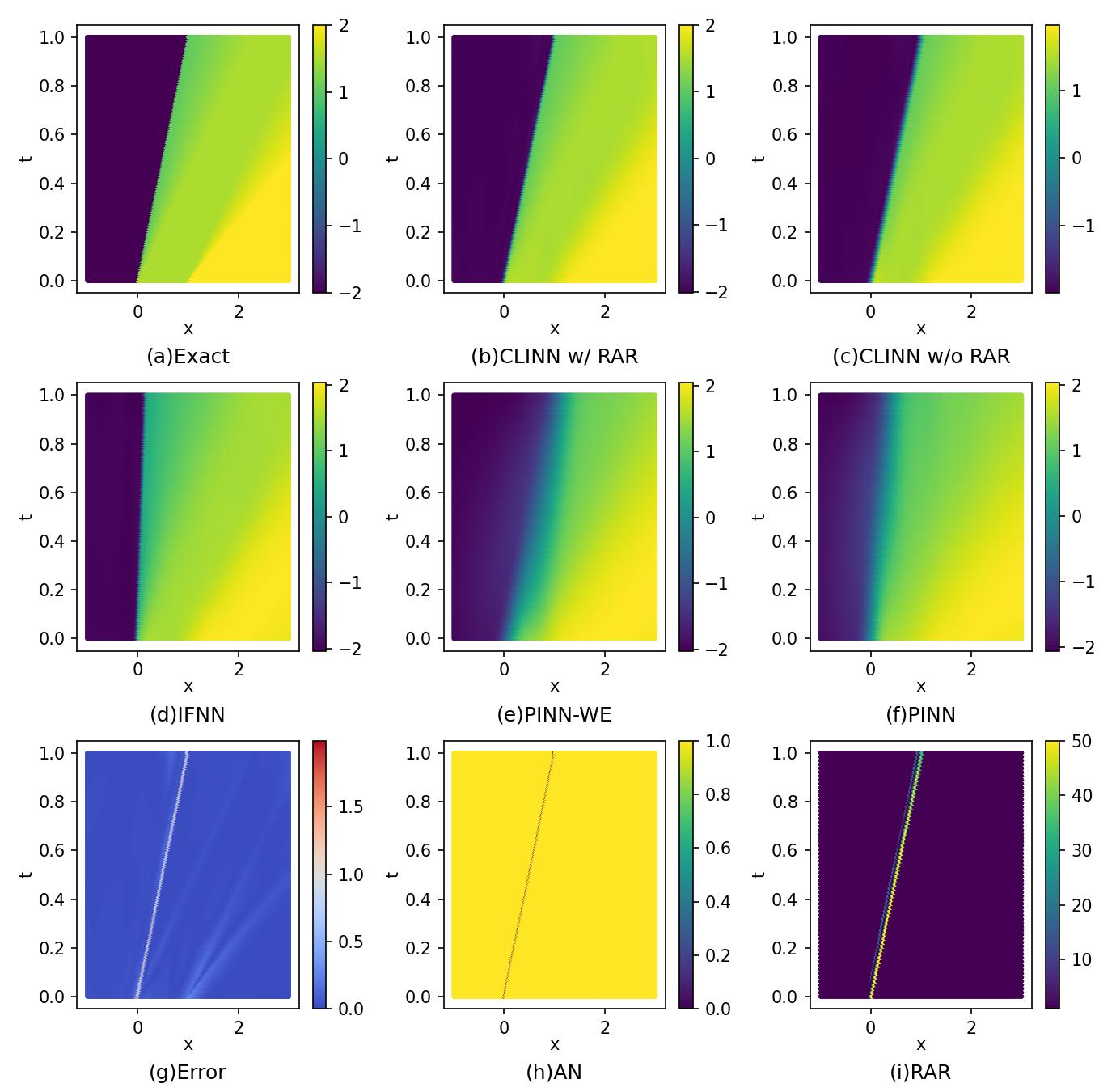}
    \includegraphics[width=\linewidth]{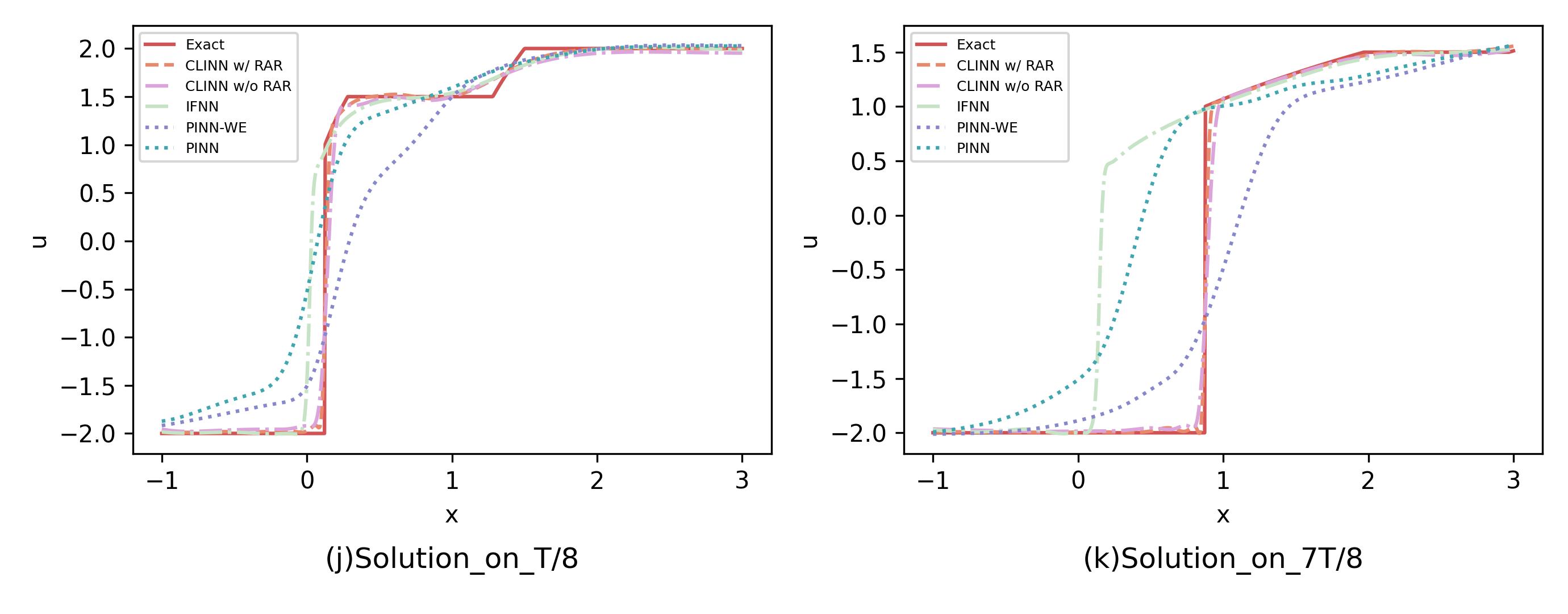}
    \caption{Performance of all methods in test case (3B). (a) Exact solution $(u_{\mathrm{true}})$; (b) Result of CLINN with RAR $(u_{\mathrm{CLINN}})$; (c) Result of CLINN without RAR; (d) Result of IFNN; (e) Result of PINN-WE; (f) Result of PINN; (g) Error of CLINN $(|u_{\mathrm{CLINN}}-u_{\mathrm{true}}|)$; (h) Discontinuity indicator result of AN; (i) Weight distribution in RAR; (j) Comparative results of all methods at $t=T/8$; (k)Comparative results of all methods at $t=7T/8$.}
    \label{case_9}
\end{figure}

Fig. \ref{case_9} compares the performance of each method on the test case (3B). Results verify that within the region $0\leq x\leq t\leq 1$, both $u=-2$ and $u=\sqrt{x/t}$ satisfy the implicit form in Eq. (\ref{solution}), but the jump condition in Eq. (\ref{RH}) uniquely determines $u=-2$ as the physically valid solution. The predictive results demonstrate that methods lacking either the implicit solution constraint or the jump condition exhibit varying degrees of distortion in this region, whereas CLINN accurately captures the solution morphology through the combined guidance. However, all methods - including CLINN - show low performance in fitting two centered rarefaction waves, as evidenced in Fig. \ref{case_9}(j).

Table \ref{Compare_MSE_3} presents the error comparison between CLINN and other advanced deep learning methods, demonstrating that CLINN achieves significantly lower errors than all three baseline approaches in both test cases. As quantified in Eq. (\ref{improve}), CLINN exhibits remarkable improvement ratios of 88.5\% and 98.0\% over PINN.

\begin{figure}[H]
    \centering
    \includegraphics[width=\linewidth]{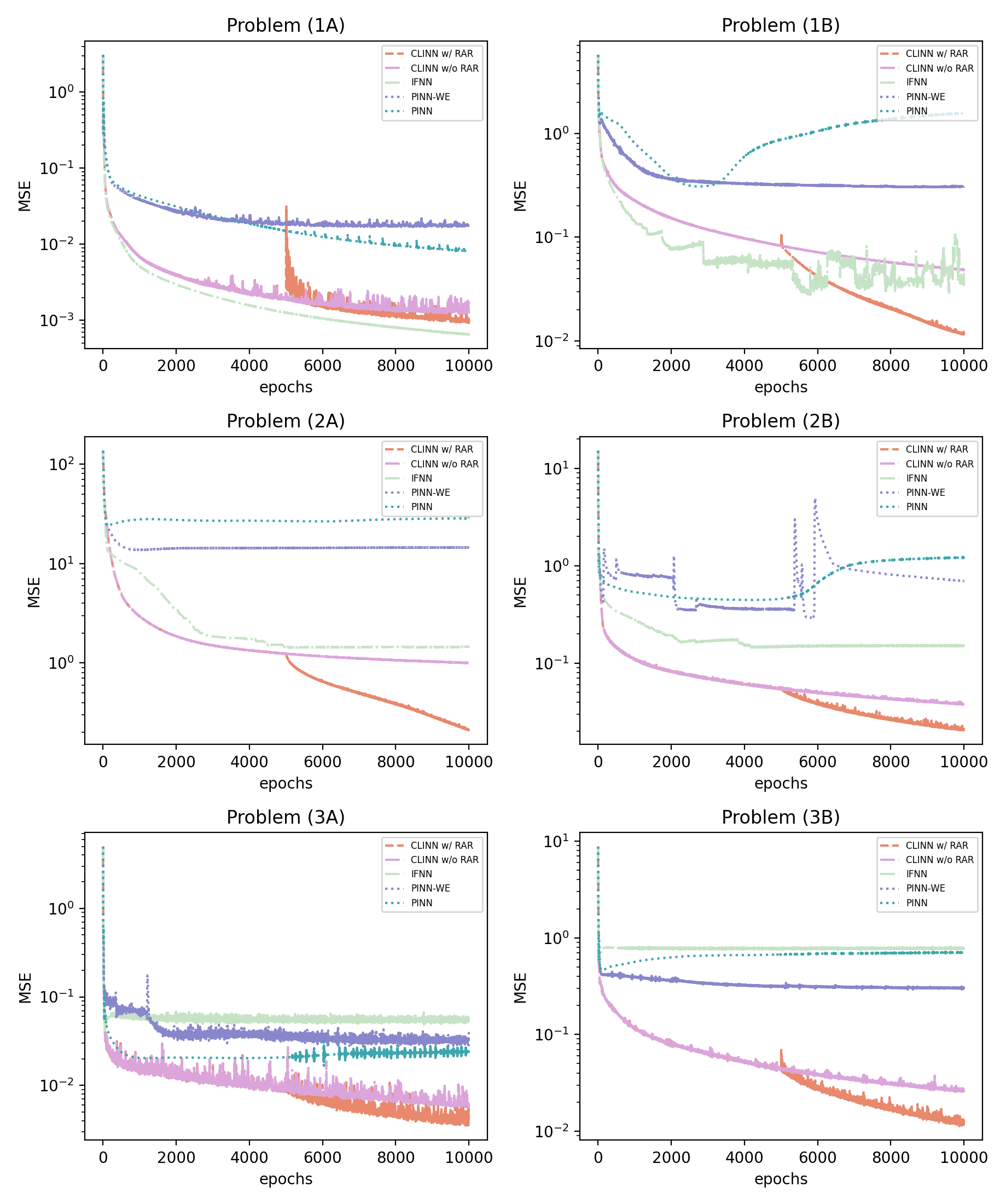}
    \caption{Comparison of MSE loss during the training process.}
    \label{loss_epoch}
\end{figure}

At the end of the research in 1D problems, we present the evolution curves of the loss versus epochs for both CLINN and other advanced deep learning methods across various test cases, as is shown in Fig. \ref{loss_epoch}. Here, the term ``loss" refers to MSE between the predicted and the exact solutions. It should be emphasized that this particular error metric is not incorporated into the total loss function for model training. Our findings reveal that CLINN achieves lower errors than other advanced deep learning methods in all test cases except for (1A). Furthermore, beyond $N_{\text{epoch},0} = 5000$ epochs, CLINN outperforms its ablation model due to the application of the weight adjustment in RAR method.

\newpage

\begin{table}[H]
\caption{Comparation among advanced deep learning methods in non-convex/non-concave problems}\label{Compare_MSE_3}
\centering
\begin{tabular}{cc|cccc|c}
\hline
Case & Model & 
MSE\_T1 &
MSE\_T2 & 
MSE\_T3 &
MSE\_T4 & 
MSE\_All
\\ \hline

\multirow{4}{*}{(3A)} 
& CLINN
& 3.36e-03
& 3.27e-03
& 3.49e-03
& 3.93e-03
& \bf{3.69e-03}
\\
& IFNN
& 1.69e-02
& 5.09e-02
& 6.84e-02
& 7.56e-02
& 5.22e-02
\\
& PINN-WE
& 2.35e-02
& 2.21e-02
& 3.10e-02
& 3.39e-02
& 2.74e-02
\\

& PINN
& 1.98e-02
& 3.61e-02
& 3.92e-02
& 3.47e-02
& 3.20e-02
\\\hline
\multirow{4}{*}{(3B)} 
& CLINN
& 1.30e-02
& 1.08e-02
& 1.02e-02
& 1.14e-02
& \bf{1.15e-02}
\\
& IFNN
& 1.90e-01
& 5.59e-01
& 9.50e-01
& 1.34e+00
& 7.62e-01
\\
& PINN-WE
& 2.64e-01
& 3.56e-01
& 3.32e-01
& 2.51e-01
& 2.98e-01
\\

& PINN
& 2.04e-01 
& 4.46e-01
& 6.52e-01 
& 9.36e-01
& 5.64e-01
\\\hline
\end{tabular}
\end{table}

\subsection{Two-Dimenstinal Problems}
\label{2D}

We consider the 2D Burgers equation, with the flux given by
\begin{equation}
    \boldsymbol f(u)=\left(\dfrac{u^2}{2}, \dfrac{u^2}{2}\right)
\end{equation}
and the initial condition is chosen to be
\begin{equation}
    u_0(x,y)=\begin{cases}
    1,&(x,y)\in(-\infty,-2]^2
    \\
    5,&(x,y)\in(-\infty,2]^2\backslash(-\infty,-2]^2\\
    -3,&(x,y)\in\mathbb{R}^2\backslash(-\infty,2]^2

\end{cases}
\label{IC-2D}
\end{equation}

To obtain the exact solution, note that when $y$ vanishes, the solution of the 1D Burgers equation is
\begin{equation}
        u_{\mathrm{1-dim}}(x,t)=5-2u_{\mathrm{2B}}(x,t)
    \label{case2-solu}
    \end{equation}
with the discontinuity curve
\begin{equation}
    x=\gamma(t)=\begin{cases}
    t+2,& t<1;\\
    -2-3t+8\sqrt{t}, &1\leq t<4; \\
    -t+6, & t>4.
    
    \end{cases}
\end{equation}

Combining the symmetry of the initial condition (Eq. (\ref{IC-2D})) with respect to the line $y=x$, it follows that the solution to the two-dimensional problem is
\begin{equation}
    u_{\mathrm{2-dim}}(x,y,t)=u_{\mathrm{1-dim}}(\max\{x,y\}, t)
\end{equation}
with the discontinuity surface
\begin{equation}
    \Phi(x,y,t):=(x-\gamma(t))(y-\gamma(t))=0,\quad x,y\leq\gamma(t).
\end{equation}
We solve the problem on $(x,y,t)\in[-4, 6]\times[-4,6]\times[0, 6] $.

\begin{figure}
    \centering
    \includegraphics[width=\linewidth]{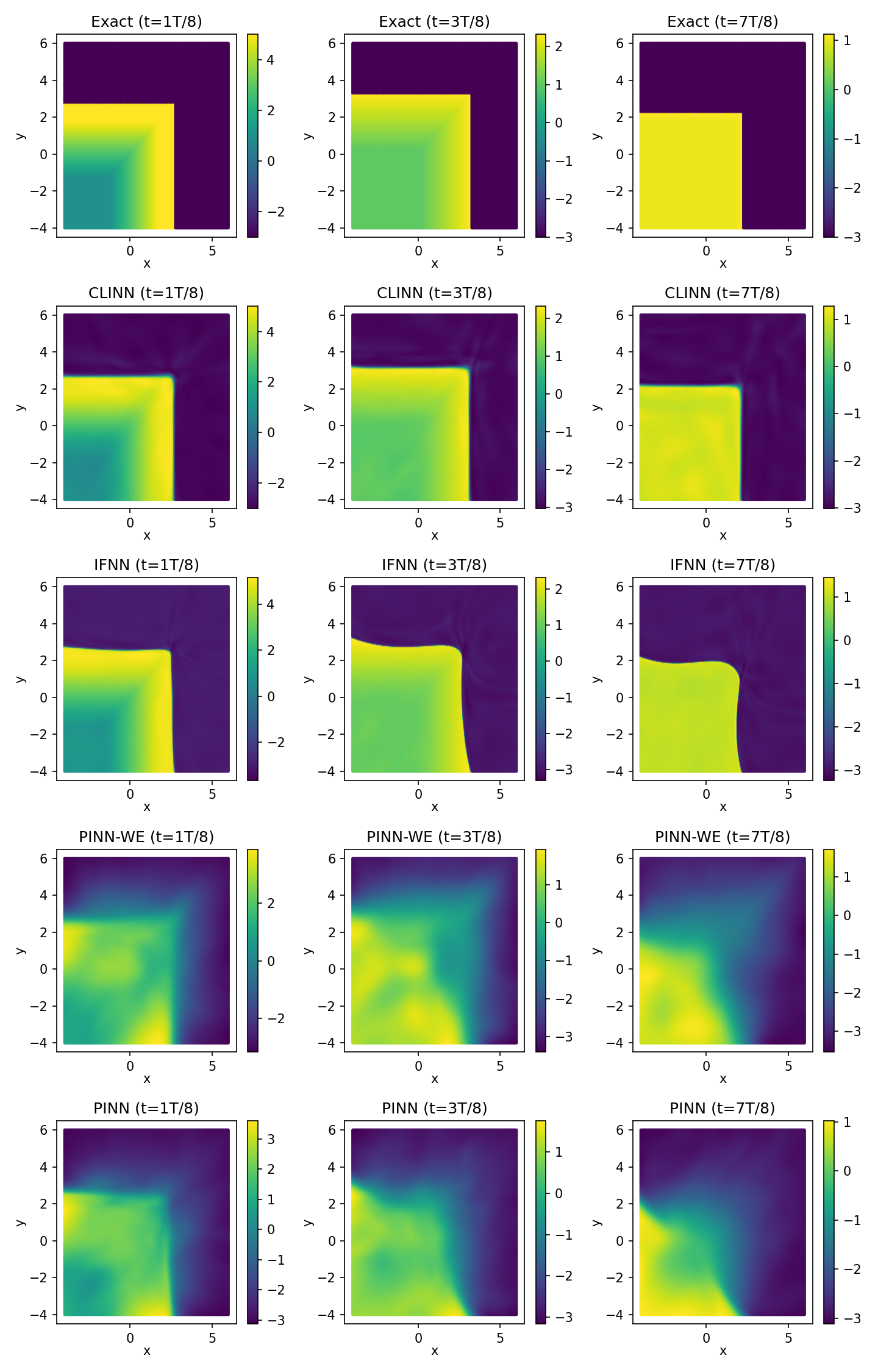}
    \vspace{-2em}
    \caption{Performance of each method in 2D Burgers equation.}
    \label{case_2d}
\end{figure}

Fig. \ref{case_2d} illustrates the results obtained by various methods. It is observed that CLINN achieves excellent agreement with the exact solution. IFNN captures the overall structure of the solution but exhibits noticeable discrepancies in representing the discontinuity. In contrast, both PINN-WE and PINN fail to accurately capture the essential features of the solution. These observations are further corroborated by the quantitative comparisons of MSE presented in Table \ref{Compare_MSE_2D}.

\begin{table}[H]
\caption{Comparation among advanced deep learning methods in 2D Burgers equations}\label{Compare_MSE_2D}
\centering
\begin{tabular}{c|cccc|c}
\hline
Model & 
MSE\_T1 &
MSE\_T2 & 
MSE\_T3 &
MSE\_T4 & 
MSE\_All
\\ \hline

CLINN
& 2.12e-01
& 1.92e-01
& 7.84e-02
& 4.51e-02
& \bf{1.08e-01}
\\

IFNN
& 9.03e-01
& 1.57e+00
& 5.64e-01
& 4.73e-01
& 8.80e-01
\\

PINN-WE
& 2.99e+00
& 1.54e+00
& 9.44e-01
& 6.96e-01
& 1.35e+00
\\

PINN

& 3.05e+00
& 1.93e+00
& 1.08e+00
& 7.54e-01
& 1.51e+00
\\\hline

\end{tabular}
\end{table}

\section{Conclusions}
\label{sec:conclusions}

We proposed CLINN, a novel neural network solver that incorporates conservation law characteristics and shock detection/correction techniques, to address the challenges of PINN in solving hyperbolic conservation laws, particularly the difficulties in preserving conservation properties and effectively approximating discontinuous solutions. The research encompassed the following two aspects:

\begin{itemize}
    \item The design of a neural network solver for scalar conservation laws. To embed conservation law features, we incorporated the boundedness constraint, implicit solution, and discontinuity jump condition into the neural network. To approximate discontinuous solutions, we dynamically adjusted the loss weights of points across the solution domain, intensifying training focus on high-error points in non-discontinuous regions.
    \item The numerical validation of test cases based on scalar conservation laws. We utilized classical benchmarks, including Burgers' equations,  Buckley-Leverett equations , and LWR models to verify the solver's effectiveness. We employed convex, concave, and non-convex/non-concave functions as numerical fluxes, along with both continuous functions and piecewise smooth functions as initial conditions, comprehensively covering fundamental solution types of the equations.
\end{itemize}

According to the result of numerical experiments in this paper, we draw the following conclusions: 

\begin{itemize}
    \item Under the combined guidance of governing equations, initial/boundary conditions, boundedness constraints, implicit solutions, and discontinuity jump conditions, CLINN demonstrates significantly enhanced approximation capability for solutions of conservation laws compared to conventional PINN.
    \item The application of the improved RAR method for dynamically adjusting loss weights across computational points yields additional improvements in the neural network's ability to capture discontinuous solutions accurately.
\end{itemize}

Building upon the current research challenges and unresolved aspects, the future work will focus on the following three directions:

\begin{itemize}
    \item Analyzing the approximation error of the neural network, to  demonstrate that CLINN achieves superior approximation performance compared to PINN theoretically.

    \item 
    Investigating under what circumstances and through what approaches to incorporate the jump conditions, to mitigate potential adverse effects on network optimization.
    \item 
    Extending the methodology proposed in this work to systems of conservation laws, to address more complex fluid dynamics problems.
\end{itemize}

\section*{Acknowledgments}
The author gratefully acknowledge the support of National Natural Science Foundation of China (No. 12571580).

\end{document}